\documentclass{tac}

\usepackage{graphicx,float}
\usepackage[utf8]{inputenc}
\usepackage{enumitem}   
\usepackage[english]{babel}
\usepackage{amsmath,amssymb,amsfonts,mathrsfs,mathtools,MnSymbol}

\usepackage{bbold}

\usepackage{quiver}

\usepackage[colorlinks=true]{hyperref}
\hypersetup{allcolors=[rgb]{0.1,0.1,0.4}}

\author{Dominik Trnka}

\thanks{Supported by RVO: 67985840 and Praemium Academi\ae of M.~Markl. Brno Ph.D.~Talent Scholarship Holder~-~Funded by the Brno City Municipality.}

\address{Department of Mathematics and Statistics, Masaryk University, Kotl\'a\v rsk\'a 2, 61137 Brno, and Institute of Mathematics, Czech Academy of Sciences, {\v Z}itn{\'a} 25, 115 67 Prague 1, Czech Republic}

\title{Category-colored operads, internal operads, and Markl $\bO$-operads}

\copyrightyear{2023}



\eaddress{trnka@math.cas.cz}

\tikzcdset{every label/.append style = {font = \normalsize}}

\newcommand{\V}{\mathscr{V}}
\newcommand{\Cat}{\it{Cat}}
\newcommand{\Endm}{\it{End}_M}
\newcommand{\Vect}{\it{Vect_{\mathbb{k}}}}
\newcommand{\NN}{\mathbb{N}}
\newcommand{\uu}{\mathbb{1}}
\newcommand{\bO}{\mathbb{O}}

\newcommand{\OO}{\mathscr{O}}
\newcommand{\SSS}{C_{\V}}
\newcommand{\T}{\mathscr{T}}
\newcommand{\sS}{\mathscr{S}}
\newcommand{\lT}{\texttt{lT}}
\newcommand{\bH}{\mathbb{H}}
\newcommand{\bS}{\mathbb{S}}

\def\colorop #1(#2;#3){
{#1}\left(
\rule{0pt}{15pt}\right.
	\hskip -3mm \begin{array}{c}
                 #3\\
		 #2
	 \end{array}
\hskip -3mm \left.
  \rule{0pt}{15pt} \right)
}
\def\coloropsq #1[#2;#3]{
{#1}\left[
\rule{0pt}{15pt}\right.
	\hskip -3mm \begin{array}{c}
                 #3\\
		 #2
	 \end{array}
\hskip -3mm \left.
  \rule{0pt}{15pt} \right]
}
\def\coloropbig #1(#2;#3){
{#1}
\begin{pmatrix}
#3 \\
#2
\end{pmatrix}
}

\let\pf\proof
\let\epf\endproof

\newtheorem{theoremA}{Theorem~A}
   \newtheorem{theoremB1}{Theorem~B1}
   \newtheorem{theoremB2}{Theorem~B2}
   \newtheorem{theoremB3}{Theorem~B3}

\begin{document}
\maketitle

\begin{abstract}
We present a Markl-style definition of operads colored by a small category. In the presence of a unit these are equivalent to substitudes of Day and Street. We show that operads colored by a category are internal algebras of a certain categorical operad of functors. We describe a groupoid-colored quadratic binary operad, whose algebras are non-unital Markl operads in the context of operadic categories. As a by-product we describe the free internal operad construction.
\end{abstract}

\section*{Introduction}\label{section:introduction}
In \cite{vdLaan:coloured} van der Laan proved that
non-symmetric Markl operads (i.e.~operads  with binary `partial
compositions' $\circ_i$, cf.~\cite[Definition 1.1.]{markl1996models}), are algebras for a colored Koszul quadratic operad with
colors the natural numbers.  A similar result for
symmetric operads was obtained by Dehling and Vallette in the
fascinating paper \cite{Dehling_Vallette:symmetric}, where 
symmetric operads were presented as algebras for a linear-quadratic curved Koszul colored operad. 
The curved Koszul theory was necessary since the authors in
loc.~cit.\ wanted to resolve, along with the binary structure
operations, also the symmetric group
actions.

Our aim is to modify van der Laan's approach to symmetric operads
and, more generally, to Markl $\bO$-operads over
an operadic category $\bO$ \cite[Definition~6.1]{Batanin_Markl:kodu2022},
considering the group (resp.~groupoid) actions as fixed parts of the structure.
The tool of choice should be operads with colors in a
groupoid that are generated, as in the van der Laan case, only by binary operations
satisfying quadratic relations.

Markl \hbox{$\bO$-operads} play an important role in homotopy theory of operad-like structures. In~\cite{Batanin_Markl:kodu2022}, various kinds of operads are viewed as algebras over the terminal Markl \hbox{$\bO$-operad}. This includes structures such as the classical operads, cyclic or modular
operads, wheeled properads, dioperads, $\frac{1}{2}$PROPs, permutads, or pre-permutads. The terminal operads for the above operad-like structures are binary quadratic. In \cite{Batanin_Markl_Obrad:models}, some of the terminal operads are resolved, providing strongly homotopy versions of the corresponding operad-like structures.

An inspiration for us was a result of
Ward~\cite{Ward:Massey}, who showed that modular operads are
algebras of a groupoid-colored quadratic Koszul operad. It
turned out that, in contrast to~\cite{Ward:Massey}, to present general Markl
$\bO$-operads as algebras of a quadratic operad, we need to work with
coloring groupoids, which have morphisms between different
objects.  For a clear exposition, we decided to study operads colored by small categories in general.
We introduce new operations $\otimes_i$ on certain functor categories, which allows us to define category-colored operads explicitly via generators and relations.
\subsection*{What is a category-colored operad?}
For a set $C$, a $C$-colored operad~$P$ consists of objects $$\colorop P(c_1\,\cdots\,c_n;c)$$ of abstract $n$-ary operations, whose inputs and output have specific types $$c_1,\ldots,c_n,c \in C.$$ The operad structure tells us how to compose two operations, provided the output of the second operation matches the type of an input of the first operation. It also tells us how to permute the inputs of an operation, and that the composition is associative and compatible with permutations. 

\begin{figure}
    \centering
    \includegraphics{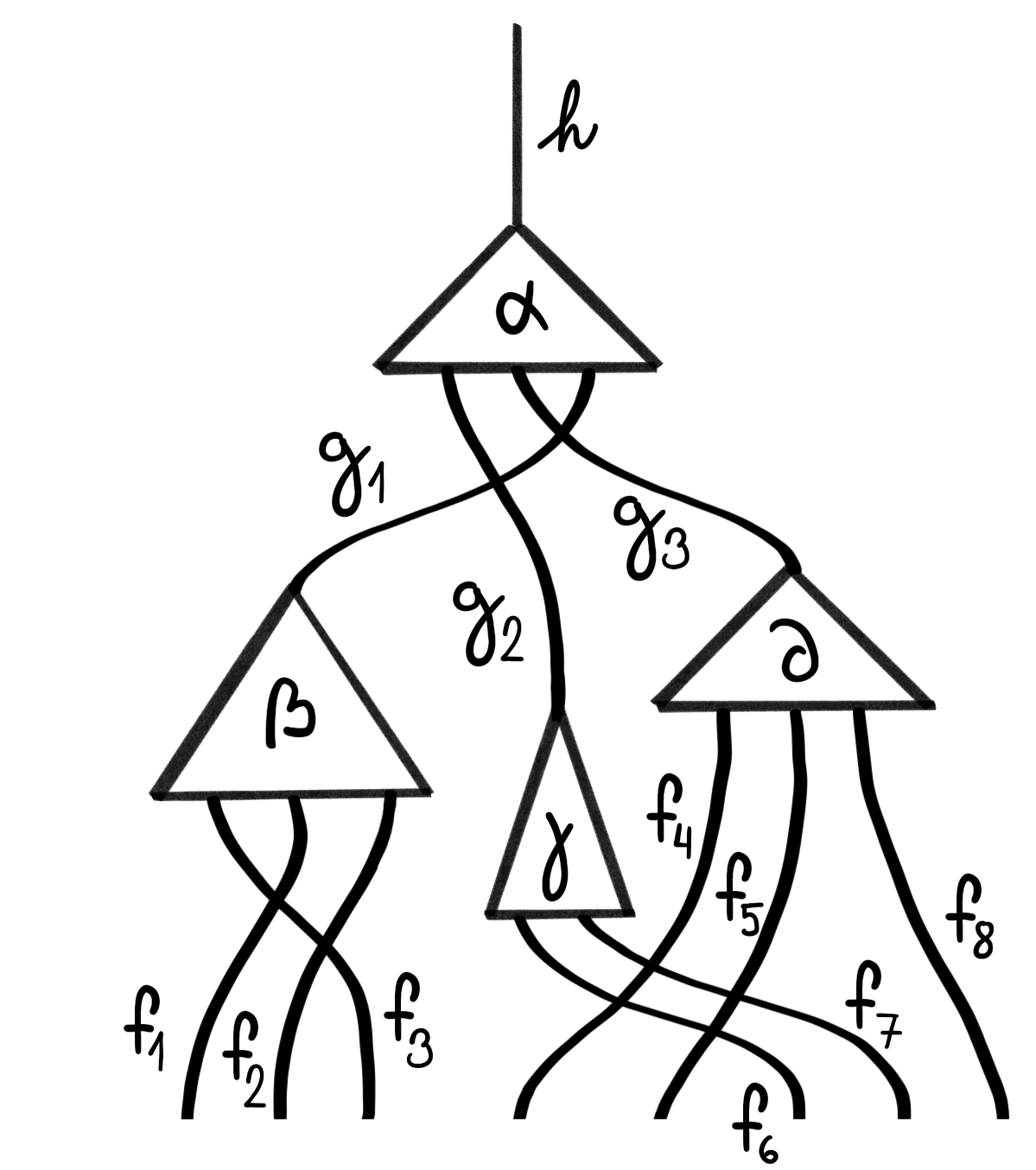}
    \caption{Composition scheme for category-colored operads.}
    \label{fig:labeled_tree}
\end{figure}

When $C$ is a small category, its morphisms act on inputs and output of the operations, possibly changing their type. Thus, two operations can be composed even if the input and output does not match, but is connected by a morphism.
Such generalisations of colored operads appeared already in \cite{Day:substitution,petersen2013operad, Ward:Massey}.
The structure of a category-colored operad is captured geometrically by labeled non-planar trees, cf.~Figure~\ref{fig:labeled_tree}. The edges are labeled by morphisms of the coloring category and vertices are labeled by operations of the operad. Such a tree determines a composition of operations and \hbox{$C$-actions} of the operad. Operads colored by categories or groupoids found applications in deformation theory \cite{DSVV}, or in homotopy theory~\cite{Batanin_White:substitudes, stoeckl2023koszul}.

One of the features of category-colored operads is the approach to unary operations of algebras. Consider the following problem: we are given an algebraic structure with unary and binary operations. We would like to describe a homotopy coherent version of this structure, but where only the binary operations are relaxed, while the unary operations remain strict. 
The solution is to encode the structure as an algebra of a category-colored operad and hide the unary operations into the coloring category. 
The problem is solved by resolving the operad in the category of colored operads. An example of such algebraic structure is a traditional operad. The unary operations are $\Sigma$-actions and the binary operations are the operad compositions $\circ_i$, cf.~Section~\ref{section:operad_for_operads}. 
Another example is a differential graded associative algebra~$A$, cf.~Example~\ref{example:dgas}. We regard the differential as a unary operation $\partial\colon A(n)\to A(n-1)$. The~coloring category $D$ has objects the integers. The Hom-spaces $D(n,n-1)$, for any integer~$n$, are generated by an element $\partial_n$, while the other Hom-spaces are trivial, which forces $\partial\partial=0$.

\subsection*{Approaches to category-colored operads:}
Classically, there are four equivalent approaches to unital operads:
\begin{enumerate}
    \item[(a)] May's definition is based on composition maps $\gamma$, which combine $n$ operations with an $n$-ary operation,
    \item[(b)] Markl's definition is based on partial compositions $\circ_i$, which plug one operation into the $i$-th input of another, 
    \item[(c)] operads are algebras of the monad of trees, and
    \item[(d)] operads are monoids in the monoidal category of $\Sigma$-modules.
\end{enumerate}
Recent approaches to operads and operad-like structures involve Feynman categories of Kaufmann and Ward~\cite{Kaufmann_Ward:Fey}, polynomial monads~\cite{Batanin_Berger:polynomial}, or operadic categories of Batanin and Markl~\cite{Batanin_Markl:duoidal}. For example, classical operads are operads over the operadic category $\it{Fin}$, the skeletal category of finite sets. 

In the category-colored case, (a)-type definition was first given by Day and Street in~\cite{Day:substitution}, under the name symmetric substitude, (c)-type definition was worked out by Ward in~\cite{Ward:Massey} for groupoids with no morphisms between different objects, and (d)-type was given by Petersen in~\cite{petersen2013operad}. 
 We complete the list by giving below Definition~\ref{definition:C_operad} of type (b). We will also address the situation when the coloring category is $\V$-enriched, cf.~Definition~\ref{definition:nu_C_operad_V}. 
  Describing category-colored operads as operads over a suitable operadic category is the subject of our current research. Category-colored operads can also be defined as monads in the bicategory of generalised species~\cite{FGHW:Species}. The formalism of operads colored by dg-categories was recently developed in~\cite{CCN:moduli}.
 
The Definition \ref{definition:C_operad} uses operations $\otimes_i$ on functor categories $${C}_{\V}(n)=\Cat\big((C^{\it{op}})^{n} \times C,\V\big), $$
where $C$ is a small category and $\V$ is a cocomplete symmetric monoidal closed category.
For functors
\[
(C^{\it{op}})^{\times n} \times C  \xrightarrow{P_n} \V,\hspace{1cm}(C^{\it{op}})^{\times m} \times C  \xrightarrow{P_m} \V,
\]
 and $1\leq i\leq n$, the functor $P_n\otimes_i P_m \in {C}_{\V}(n+m-1)$ is a colimit coequalizing $C$-actions on the $i$-th input of $P_n$ and output of $P_m$.
The $C$-operad is then a collection of functors $P_n$ with natural transformations 
\begin{equation*}
       \begin{tikzcd}
P_n\otimes_i P_m & P_{m+n-1},
	\arrow["\circ_i",from=1-1, to=1-2]
\end{tikzcd}\end{equation*}
satisfying classically-looking associativity and equivariance axioms.
Let us describe two main results of this article.

\subsection*{Result I}
We give a conceptual explanation of the operations $\otimes_i$. In Section~\ref{section:SSS} (Proposition~\ref{proposition:C_V_is_an_operad}) we prove that the operations $\otimes_i$ equip the collection of categories $C_{\V}(n)$ with an operad structure. However, the associativity of $\otimes_i$ holds only up to canonical isomorphisms. More precisely, $C_{\V}$ is an example of a (non-strict) categorical operad, that is, an operad with values in the monoidal category of categories, where the operad axioms hold only up to coherent natural isomorphisms. Categorical operads will be defined in Section~\ref{section:internal_operads} as pseudo-algebras of a 2-monad of trees. 

For a 2-monad $T$, an internal algebra of a pseudo-$T$-algebra $\OO$ is a lax morphism from the terminal pseudo-$T$-algebra to $\OO$. A nice example of internal algebras are monoids in monoidal categories:
a monoidal category $M$ is a pseudo-$\mathbb{M}$-algebra of the 2-monad $\mathbb{M}$ on the 2-category of categories, which is induced from the `free monoid' monad in sets. A~monoid $m$ in $M$ is then an internal algebra in the pseudo-$\mathbb{M}$-algebra $M$. Note that monoidal categories provide just enough coherent data for a definition of a monoid. From this point of view, categorical operads provide the right coherent data for a definition of operads. Since monoidal categories are categorical operads concentrated in arity 1, our philosophy is sound with the classical definition of operads in monoidal categories. 
Our first result is:
\begin{center}
   \textsc{Theorem A} (\ref{proposition:iso_of_int_operads_and_C_operads}). \textit{$C$-operads are internal operads in the categorical operad~$(C_{\V}(n),\otimes_i)$.} 
\end{center}
Further, we describe free internal operads and, using Theorem A, we give an explicit free $C$-operad construction, which will be used in Result II. However, the free construction can be read off the monoidal or monadic definitions of \cite{petersen2013operad,Ward:Massey} and hence the Result II can be understood independently of Result I.
\subsection*{Result II}
Our second result is based on ideas of~\cite{Batanin_Markl_Obrad:models}, that express operads as algebras over a groupoid-colored operad (we will use the term hyperoperad), and also on the above mentioned result of Ward~\cite{Ward:Massey}. The coloring groupoid incorporates the symmetric action and equivariance axiom of operads to the level of collections. Hence, we do not regard the symmetric actions as a piece of operadic structure, which simplifies the process of resolving the 
 hyperoperad significantly. This should be put in contrast with~\cite{Dehling_Vallette:symmetric}, which treats the symmetric actions as unary operations, hence resolving them as well. For more background and history on resolving hyperoperads, we refer the reader to the introduction of \cite{Batanin_Markl_Obrad:models}.

The groupoid of colors for the hyperoperad $\bH$, whose algebras are classical non-unital symmetric Markl operads, will be denoted by $\Sigma$. It has natural numbers as objects and permutations as morphisms.
The idea is to generate the hyperoperad $\bH$ by abstract symbols $$*_i \in \colorop \bH(n\hspace{.5cm}m;n+m-1),$$ which represent the
operad structure maps
\begin{equation}\label{eq:circ_i_intro}
       \begin{tikzcd}
O(n)\otimes O(m) & O({m+n-1}).
	\arrow["\circ_i",from=1-1, to=1-2]
\end{tikzcd}\end{equation} 
We equip the generating collection with free actions of $\Sigma$ on the inputs and the output, obtaining formal elements $\colorop (*_i,\sigma\,\tau;\hspace{15pt}\pi)$, which represent the composites 
\[
\begin{tikzcd}
	 O(n)\otimes O(m) & O(n+m-1)\\ O(n)\otimes O(m)  & O(n+m-1)
	\arrow["(-)\sigma \otimes (-)\tau", from=2-1, to=1-1]
 \arrow["\circ_i", from=1-1, to=1-2]
 \arrow["(-)\pi", from=1-2, to=2-2]
 \arrow[dashed, from=2-1, to=2-2]
\end{tikzcd}
\]
of symmetric actions and compositions $\circ_i$ of a classical operad $O$.
On the level of collections we make the symbols $*_i$ equivariant by identifying
$$\colorop (*_i,\sigma\,\tau;\hspace{15pt}1)\sim \colorop(*_{\sigma(i)},1\hspace{10pt}1;\hspace{28pt}\sigma\circ_i\tau).$$
It remains to generate the free operad from the constructed collection and divide by the operadic ideal generated by associators, which ensures associativity of the symbols $*_i$. The result~is:
\begin{center}\textsc{
Theorem B1} (\ref{theorem:operads_as_algebras}).
\textit{(Symmetric non-unital) Markl operads are algebras of the binary quadratic operad $\mathbb{H}$ with colors $\Sigma$.}
\end{center}
We regard units for operads as an extra structure. To encode units $\mathbb{k}\xrightarrow{u}O(1)$ of a unital operad $O$, one uses constants ($0$-ary operations)
$$u\in \colorop \bH(\emptyset;1).$$ Notice however, that the unit axioms
$u\circ_1\alpha=\alpha$ and $\alpha \circ_i u=\alpha$ are not quadratic.
We show, that Theorem B1 extends to the framework of operads in operadic categories and $C$-colored operads by proving:
\begin{center}\textsc{
Theorem B2} (\ref{theorem:O-operads_as_algebras}).
\textit{(Non-unital) Markl $\bO$-operads in an operadic category $\bO$ are algebras of the binary quadratic operad $\mathbb{H}_{\bO}$ with colors $\bO_{\it{iso}}^{\it{op}}$, the category of isomorphisms~of~$\bO$.}
\end{center}
\begin{center}\textsc{
Theorem B3} (\ref{theorem:C-operads_as_algebras}).
\textit{For a fixed small category $C$, (non-unital) $C$-operads are algebras of the binary quadratic operad $\mathbb{H}_{C}$ with colors $Bq^{\Sigma}(C)$.}
\end{center}
 
In Section \ref{section:algebras} we observe that algebras of operads with values in a functor category~$\V^C$ can be expressed as algebras of $C$-operads in $\V$.
In our future work we shall address Koszulness of the constructed hyperoperads $\mathbb{H}_{\bO}$ and their resolutions.
\subsection*{Plan of the paper:}

The first two sections lead to a Markl-style definition of $C$-operads and  comparison with symmetric substitudes. The rest of the paper is divided into two independent parts.
The first part consists of Sections~\ref{section:internal_operads}-\ref{section:free}, which are devoted to Result~I.
The second part consists of Sections~\ref{section:algebras} and \ref{section:operad_for_operads}, where we construct the hyperoperads of Result~II. We believe that the reader can understand the Result~II without the details of the free construction from Section~\ref{section:free}, or look at Example~\ref{example:free_operad}, if necessary.

In more detail, Section~\ref{section:internal_operads} recalls the theory of polynomial monads and internal algebras. We define internal operads in a categorical operad and describe them explicitly.
    In Section~\ref{section:SSS} we prove Theorem A.
    Section~\ref{section:free} gives the free internal operad construction. Using the Theorem~A, we describe the free $C$-operad.
    Section~\ref{section:algebras} defines algebras of \hbox{$C$-operads} and adapts some standard notions to the category-colored setting. In the second part of Section~\ref{section:algebras} we compare algebras of classical operads with values in functor categories $\V^C$ with algebras of $C$-colored operads in $\V$. An application of this comparison is a description of differential graded associative algebras as algebras over a category-colored operad with values in $\Vect$.
    The last section is devoted to proving Theorems~B1-B3.
\subsection*{Acknowledgement}
I would like to express gratitude to my supervisor Martin Markl for his initial idea and constant support, and to my colleagues Miloslav Štěpán and Michael Batanin for their helpful comments. I also thank to an anonymous referee for valuable comments and suggestions.
\subsection*{Conventions}
The identity morphism of an object $x$ will be denoted by $\mathbb{1}_x$ or only by $\mathbb{1}$, leaving the object implicit. The set of natural numbers (including 0) will be denoted by~$\NN$.
If not stated otherwise, $C$ denotes a small category and $\V$ a cocomplete symmetric monoidal closed category. However, some constructions and results may work with weaker assumptions on $\V$.
\section{The $\otimes_i$ products}\label{section:notation}
We establish some notation related to functors
\[\begin{tikzcd}
\underbrace{C^{\it{op}} \times \cdots \times C^{\it{op}}}_{n\textit{-times}} \times C  &  \V, n\geq 0
	\arrow["X_n",from=1-1, to=1-2]
\end{tikzcd}\]
and their natural transformations, where $C$ is a small category and $\V$ a cocomplete symmetric monoidal closed category. We denote by $\Cat$ the 2-category of categories.
\begin{definition}\label{definition:C_V}
We denote the functor category $\Cat\big((C^{\it{op}})^{n} \times C,\V\big)$ by ${C}_{\V}(n)$. A~collection of functors 
$\{X_n\in {C}_{\V}(n)\}_{n\in\NN}
$ will be called a non-symmetric $C$-collection in $\V$. 
\end{definition}
The objects of $(C^{\it{op}})^{n} \times C$ are denoted as schemes 
$\colorop (c_1 \,\cdots\, c_n;c)$
and the maps are
\[\begin{tikzcd}
\colorop (f_1 \,\cdots \,f_n;f)\colon \colorop (d_1  \,\cdots \,d_n;d)
& \colorop (c_1 \,\cdots\, c_n;c).
	\arrow[from=1-1, to=1-2]
\end{tikzcd}\]
Fixing an index $i\in \{1,\ldots,n\}$ and an object $c \in C$, we
partially evaluate a functor \hbox{$X\in {C}_{\V}(n)$} in the $i$-th input to obtain a functor
\[\begin{tikzcd}
X^i_c = \colorop X(\bullet \,\cdots\, c \,\cdots\,
\bullet;\bullet)\colon \underbrace{C^{\it{op}} \times \cdots \times  C^{\it{op}}}_{(n-1)\textit{-times}} \times C  &  \V.
	\arrow[from=1-1, to=1-2]
\end{tikzcd}\]
Any $f \colon d \rightarrow c$ in $C$ then induces a natural transformation 
\[\begin{tikzcd}
X^i_c   &  X^i_d.
	\arrow["X^i_f",from=1-1, to=1-2]
\end{tikzcd}\]
We thus have a functor $X^i_{\bullet}$ from $C^{\it{op}} $ into the functor category
\[\begin{tikzcd}
C^{\it{op}}   &  \Cat(\underbrace{C^{\it{op}} \times \cdots \times  C^{\it{op}}}_{(n-1)\textit{-times}}\times C, \V).
	\arrow["X^i_{\bullet}",from=1-1, to=1-2]
\end{tikzcd}\]
Similarly, the partial evaluation $\prescript{c}{}{X}:= \colorop X(\bullet\, \cdots\,\bullet;c)$ gives a functor

\[\begin{tikzcd}
C   & \Cat(\underbrace{C^{\it{op}} \times \cdots \times  C^{\it{op}}}_{n\textit{-times}}, \V).
	\arrow["\prescript{\bullet}{}{X}",from=1-1, to=1-2]
\end{tikzcd}\]
For two functors  $X\in C_{\V}(n), Y \in C_{\V}(m),$ define a functor $$X^i_c\otimes \prescript{d}{}{Y} \in {C}_{\V}(n+m-1)$$ by the formula $$\colorop {(X^i_c\otimes \prescript{d}{}{Y})} (x_1\,\cdots\,y_1\,\cdots\, y_m\, \cdots \,x_{n};x) :=\colorop X (x_1\,\cdots \,x_{i-1}\hspace{5pt} c \hspace{5pt} x_{i+1}\,\cdots\, x_{m};x)\otimes \colorop Y(y_1\,\cdots\, y_{m};d).$$
This induces a functor 
\begin{equation}\label{eq:Xi_otimes_Y}
\begin{tikzcd}
C^{\it{op}}\times C   && C_{\V}(n+m-1).
	\arrow["X^i_{\bullet}\otimes \prescript{\bullet}{}{Y}",from=1-1, to=1-3]
\end{tikzcd}
\end{equation}

\begin{lemma}\label{lemma:associativity_of_X^i_otimes_Y}
    Let
$X \in \SSS(n)$, 
$Y \in \SSS(m)$, 
$Z \in \SSS(k)$,
and $1\leq i< j\leq n$.
There are canonical natural isomorphisms
    \begin{equation}\label{eq:random2}
    (X^j_{\bullet}\otimes \prescript{\bullet}{}{Y})^i_{\bullet}\otimes \prescript{\bullet}{}{Z} \cong
    \begin{cases*}
    (X^i_{\bullet}\otimes \prescript{\bullet}{}{Z})^{j+k-1}_{\bullet}\otimes \prescript{\bullet}{}{Y}
      & if $1 \leq i< j\leq n$, \\
      X^j_{\bullet}\otimes \prescript{\bullet}{}{(Y^{i-j+1}_{\bullet}\otimes \prescript{\bullet}{}{Z})}
      & if $j \leq i<j+m$, \\
       (X^j_{\bullet}\otimes \prescript{\bullet}{}{Z})^{i-m+1}_{\bullet}\otimes \prescript{\bullet}{}{Y}
      & if $j+m \leq i\leq n+m-1$. 
    \end{cases*}
  \end{equation}
\end{lemma}
\pf
    It is enough to evaluate both sides and compare the values in $\V$. The isomorphisms~(\ref{eq:random2}) are induced by the symmetry and associativity isomorphisms of $\V$. 
\epf
Let us recall the definition of cowedges and coends:
\begin{definition}
    Let $F\colon C^{\it{op}}\times C \rightarrow D$ be a functor and $\gamma_c\colon F(c,c)\rightarrow t$ a family of maps in $D$ with a common target $t \in D$, indexed by objects $c\in C$.
The family $\gamma$ is a cowedge from the functor $F$, if the diagram 
\[\begin{tikzcd}
F(d,c) & F(d,d)
\\
F(c,c) & t
\arrow["{F(\uu,f)}",from=1-1, to=1-2]
\arrow["{F(f,\uu)}"',from=1-1, to=2-1]
\arrow["\gamma_c",from=1-2, to=2-2]
\arrow["\gamma_d",from=2-1, to=2-2]
\end{tikzcd}\]
commutes for any $f\colon c\rightarrow d$ in $C$.
\end{definition}
\begin{definition}
    The coend of a functor $F\colon C^{\it{op}}\times C \rightarrow D$ is the initial cowedge from~$F$:
 \[\begin{tikzcd}
F(c,c)  &  \displaystyle\int^{c\in C} F(c,c).
	\arrow["\iota_c",from=1-1, to=1-2]
\end{tikzcd}\]
\end{definition}
If $D$ has colimits, the coend is the colimit of the diagram
\begin{equation}\label{diagram:coend_as_colimit}
\begin{tikzcd}
	{\displaystyle\bigoplus_{f\colon c' \rightarrow c''} F(c'',c') } &&& {\displaystyle\bigoplus\limits_{c\in C} F(c,c).}
	\arrow["{\displaystyle\bigoplus_{f\colon c' \rightarrow c''} F(f,\uu)}", curve={height=-18pt}, from=1-1, to=1-4]
	\arrow["{\displaystyle\bigoplus_{f\colon c' \rightarrow c''} F(\uu,f)}"', curve={height=18pt}, from=1-1, to=1-4]
\end{tikzcd} 
\end{equation}
The symbol $\bigoplus$ stands for the categorical coproduct.
\begin{definition}\label{definition:odot_i} 
Let $X\in C_{\V}(n), Y \in C_{\V}(m),$ and $i\in \{1,\ldots,n\}$.
We define the functor $$X\otimes_i Y\in C_{\V}(n+m-1)$$ as the coend
$$
X\otimes_i Y:=\displaystyle\int^{c\in C} X^i_{c}\otimes \prescript{c}{}{Y}.
$$
\end{definition}
Any two natural transformations $\alpha\colon X'\to X''$ and $\beta\colon Y'\to Y''$ induce a natural transformation $\alpha\otimes_i\beta\colon X'\otimes_i Y'\to X''\otimes_i Y''$.  In fact, $-\otimes_i-$ is a functor 
\[\begin{tikzcd}
 \SSS(n)\times \SSS(m)  & \SSS(n+m-1).
	\arrow["\otimes_i",from=1-1, to=1-2]
\end{tikzcd}\]
The operations $\otimes_i$ are partial composition analogs of the plethysm operation, which is the monoidal product of the monoidal category of $C\mathbb{S}$-modules in $\V$ of~\cite{petersen2013operad}. 
By~definition, a~natural transformation 
 \[\begin{tikzcd}
X\otimes_i Y  &  Z
	\arrow["f",from=1-1, to=1-2]
\end{tikzcd}\]
is a collection of natural transformations 
\[\begin{tikzcd}
X^i_{c}\otimes \prescript{c}{}{Y}  &  Z,
	\arrow["f_c",from=1-1, to=1-2]
\end{tikzcd}\]
for all $c\in C$, which form a cowedge.

\begin{lemma}\label{lemma:otimes_associativity}
The operations $\otimes_i$ are associative, i.e.~for $X$, $Y$, $Z$, $i$, and $j$ as in Lemma~\ref{lemma:associativity_of_X^i_otimes_Y},
there are canonical natural isomorphisms
    \begin{equation}\label{eq:random1}
    (X\otimes_j Y)\otimes_{i} Z \cong
    \begin{cases*}
    (X\otimes_{i}Z)\otimes_{j+k-1}Y
      & if $1 \leq i< j\leq n$, \\
      X\otimes_{j}(Y\otimes_{i-j+1}Z)
      & if $j \leq i<j+m$, \\
       (X\otimes_{j}Z)\otimes_{i-m+1} Y
      & if $j+m \leq i\leq n+m-1$. 
    \end{cases*}
  \end{equation}
\end{lemma}
\pf
    The proof follows from Lemma \ref{lemma:associativity_of_X^i_otimes_Y} and formal manipulations with coends.
\epf
\begin{example}\label{example:explicit_odot}
For $\V=\Vect$, the category of vector spaces over a field ${\mathbb{k}}$, the $\otimes_i$ operations can be described as follows.

Let $X \in \SSS(n)$, $Y \in \SSS(m)$. Pick two vectors 
\begin{center}
    \begin{tabular}{c c}
      $ \alpha \in \colorop X(x_1 \,\cdots \,c \,\cdots\, x_n;x),$  &  
      $ \beta \in \colorop Y(y_1\, \cdots\, y_m;d),$
    \end{tabular}
\end{center}
 and a map $f \colon d \rightarrow c$ of $C$. The map $f$ acts on the $i$-th input of $X$ by 
\[\begin{tikzcd}
X\colorop (\uu_{x_1}\,\cdots \,f\,\cdots \,\uu_{x_n};\uu_{x})\colon X\colorop (x_1 \,\cdots\, c\, \cdots\, x_n;x)  & X\colorop (x_1\, \cdots\, d\, \cdots\, x_n;x).
	\arrow[from=1-1, to=1-2]
\end{tikzcd}\]
Similarly, $f$ acts on the output of $Y$ by 
\[\begin{tikzcd}
Y\colorop (\uu_{y_1}\,\cdots\, \uu_{y_m};f)\colon \colorop Y(y_1\, \cdots\, y_m;d)  & \colorop Y(y_1\, \cdots\, y_m;c).
	\arrow[from=1-1, to=1-2]
\end{tikzcd}\]
 Denote by $\alpha \cdot f$ and $f \cdot \beta$ the values 
    $$X\colorop (\uu\,\cdots\, f\,\cdots\, \uu;\uu)(\alpha) \in \colorop X(x_1\, \cdots \,d\, \cdots\, x_n;x) \text{ and }Y\colorop (\uu\,\cdots\, \uu;f)(\beta) \in \colorop Y(y_1 \,\cdots \, y_m;c).$$
Then each space $$\colorop X\otimes_i Y (x_1\,\cdots\, y_1\,\cdots \, y_m \,\cdots\, x_n;x)$$
is the direct sum of spaces 
$$\displaystyle\bigoplus_{c\in C}\colorop X(x_1\,\cdots\,c\,\cdots\, x_n;x)\otimes \colorop Y(y_1\,\cdots\, y_m;c),$$
divided by the ideal generated by elements $(\alpha\cdot f\otimes \beta - \alpha\otimes f\cdot \beta)$ for all $\alpha$, $\beta$ and~$f$.
\end{example}

\section{Category-colored operads}\label{section:C-operads}
\begin{figure}[H]
    \centering
    \includegraphics{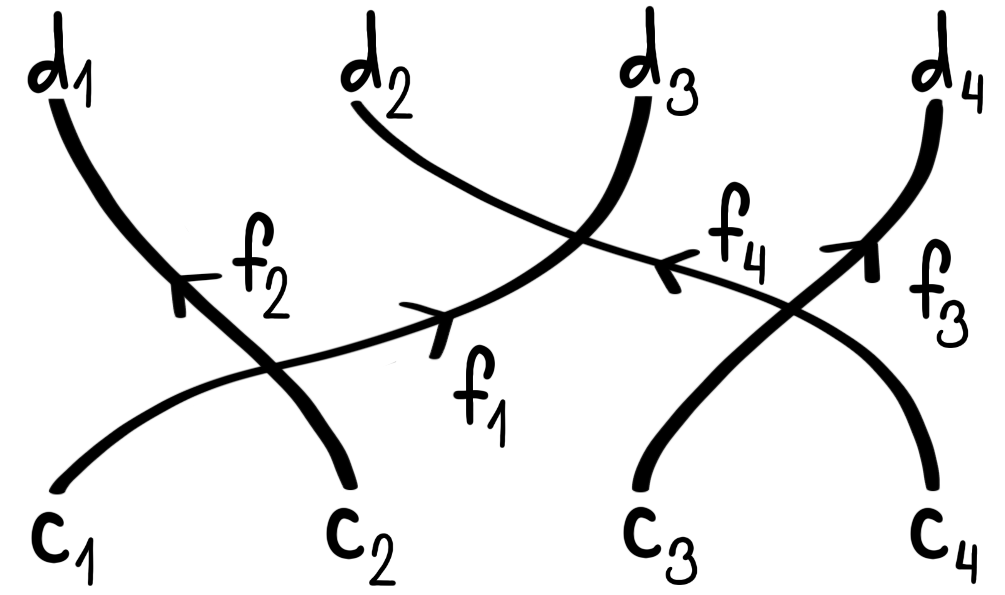}
    \caption{A map of $\bS C$.}
    \label{fig:morphism}
\end{figure}

\begin{figure}[H]
    \centering
     \includegraphics{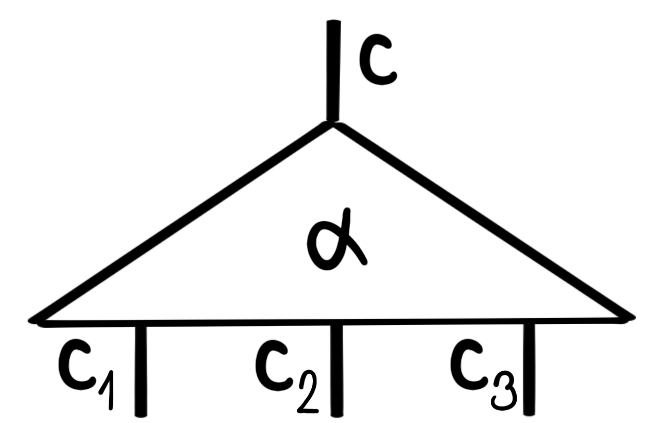}
    \caption{Abstract $n$-ary operation.}
    \label{fig:operation}
\end{figure}

Before we define operads colored by a small category $C$, we introduce their underlying (symmetric) $C$-collections.
Let $\bS C$ be the free symmetric monoidal category on $C$. The objects of $\bS C$ are ordered lists $\underline{c}=(c_1\,\cdots\,c_n)$ of objects of $C$ and maps $\underline{c}\rightarrow \underline{d}$ are tuples~$(\sigma\,f_1\,\cdots\,f_n)$, where $\sigma \in \Sigma_n$ and each $f_i\colon c_i \longrightarrow d_{\sigma(i)}$ is a map in $C$, cf.~Figure~\ref{fig:morphism}.
\begin{definition}\label{lemma:C-coll2}
    A $C$-collection in $\V$ is a functor 
\begin{equation}\label{eq:C-collection}
    \begin{tikzcd}
(\bS C)^{\it{op}} \times C 
  &  \V.
	\arrow["P",from=1-1, to=1-2]
\end{tikzcd}
\end{equation}
  Morphisms of $C$-collections are their natural transformations.
\end{definition}
   The value $\colorop P(c_1 \,\cdots\, c_n;c)$ of a $C$-collection $P$ is interpreted as a collection of abstract $n$-ary operations whose inputs are of types $c_1 ,\ldots, c_n$ and the output is of type $c$. Such an operation is depicted in Figure~\ref{fig:operation}.

\begin{figure}[H]
    \centering
    \includegraphics[]{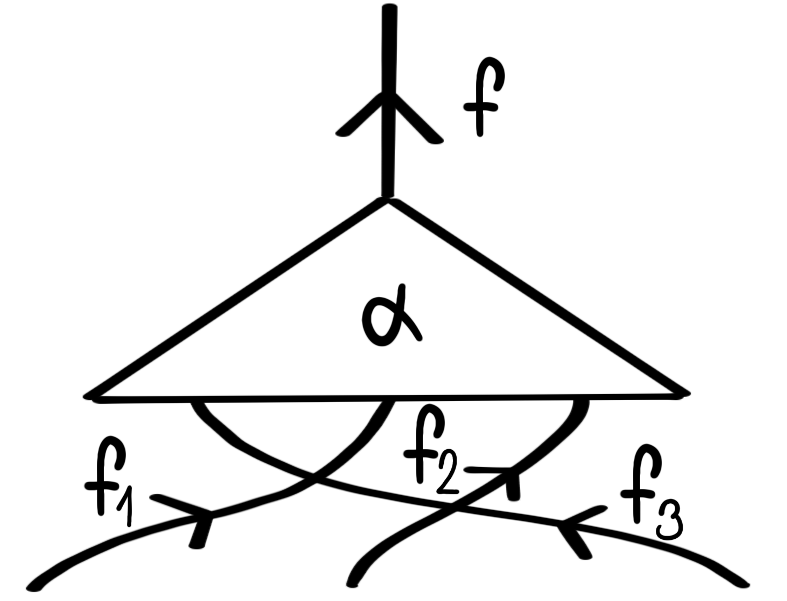}
    \caption{Acting on an operation.}
    \label{fig:acting_on_operation}
\end{figure}

The morphisms of $(\bS C)^{\it{op}} \times C$ are
\[\begin{tikzcd}
\colorop (\sigma\, f_1 \,\cdots \,f_n;\hspace{8pt}f)\colon \colorop (d_1  \,\cdots \,d_n;d)
& \colorop (c_1 \,\cdots\, c_n;c)
	\arrow[from=1-1, to=1-2]
\end{tikzcd}\]
and the value $\colorop P(\sigma\, f_1 \,\cdots \,f_n;\hspace{8pt}f)$ is interpreted as a permutation of inputs, together with an~action of maps of $C$ on inputs and output of operations. The actions thus change the type of an~operation, cf.~Figure~\ref{fig:acting_on_operation}. 
Any morphism $(\sigma\, f_1\,\cdots\,f_n)$ of $\bS C$ can be decomposed as
\[
\begin{tikzcd}
	& {(d_{\sigma(1)}\,\cdots\,d_{\sigma(n)})} \\
	{(c_1\,\cdots\,c_n)} && {(d_1\,\cdots\,d_n)} .\\
	& {(c_{\sigma^{-1}(1)}\,\cdots\,c_{\sigma^{-1}(n)})}
	\arrow["{(\sigma\, f_1\,\cdots\,f_n)}", from=2-1, to=2-3]
	\arrow["{(\sigma\, \uu\,\cdots\,\uu)}"', from=2-1, to=3-2]
	\arrow["{(\uu\, f_1\,\cdots\,f_n)}", from=2-1, to=1-2]
	\arrow["{(\sigma\, \uu\,\cdots\,\uu)}", from=1-2, to=2-3]
	\arrow["{(\uu\,f_{\sigma^{-1}(1)}\,\cdots\,f_{\sigma^{-1}(n)})}"', from=3-2, to=2-3]
\end{tikzcd}
\]
   A $C$-collection $P$ is therefore a non-symmetric $C$-collection of functors $P_n\in C_{\V}(n)$, for every $n\geq 0$, equipped with natural transformations 
\[\begin{tikzcd}
\colorop P_n(c_1 \,\cdots\,
c_n;c)  &  \colorop P_n(c_{\sigma(1)} \,\cdots\,
c_{\sigma(n)};c),
	\arrow["(-)\sigma",from=1-1, to=1-2]
\end{tikzcd}\]
for every permutation $\sigma\in\Sigma_n$, such that $((-)\tau)\sigma=(-)\tau\sigma$ and $(-)\uu_n=\uu_{C_{\V}(n)}$. We~will sometimes omit the subscripts of the functors $P_n$, when they are clear from the context. 
   For any $\sigma \in \Sigma_n$, $\tau \in 
    \Sigma_m$ and $c\in C$, the actions 
    $(-)\sigma \text{ and } (-)\tau$   
  on functors $X\in C_{\V}(n)$, \hbox{$Y \in C_{\V}(m)$}, induce a natural transformation 
   \begin{equation*}
       \begin{tikzcd}
X^{\sigma(i)}_{c}\otimes \prescript{c}{}{Y} &&&  X^i_{c}\otimes \prescript{c}{}{Y},
	\arrow["(-)\sigma\otimes (-)\tau",from=1-1, to=1-4]
\end{tikzcd}
   \end{equation*}
   which extends to colimits
   \begin{equation*}
       \begin{tikzcd}
X\otimes_{\sigma(i)} Y &&&  X\otimes_i Y.
	\arrow["(-)\sigma\otimes (-)\tau",from=1-1, to=1-4]
\end{tikzcd}
   \end{equation*}
In the following definition, the symbol `$\cdot$' denotes ordinary composition of morphisms.
\begin{definition}\label{definition:C_operad}  
A (non-unital) $C$-operad in $\V$ is a $C$-collection $P$ in $\V$,
 together with natural transformations
 \begin{equation}\label{eq:circ_i}
       \begin{tikzcd}
P_n\otimes_i P_m & P_{m+n-1},
	\arrow["\circ_i",from=1-1, to=1-2]
\end{tikzcd}\end{equation}
for any $n\geq1,m\geq0,i\in\{1,\ldots,n\}$, called the composition maps, which are required to satisfy the following two axioms.
\begin{itemize}
    \item Equivariance:
Let $\sigma\circ_i\tau$ be the permutation given by inserting the permutation $\tau$ at the $i$-th place in $\sigma$ (see \cite[II.1.3(1.8)]{markl2002operads}). For every $\sigma\in \Sigma_n$ and $\tau\in \Sigma_m$,
\begin{equation}\label{eq:equivariance} \circ_i\cdot ((-)\sigma\otimes(-)\tau )=(-)(\sigma\circ_i\tau)\cdot \circ_{\sigma(i)}.
 \end{equation}
\item Associativity: For any $1\leq j\leq n, 0\leq m$ and $0\leq k$, 
\begin{equation}\label{eq:associativity}
    \circ_i\cdot (\circ_j\otimes_i \uu)=
    \begin{cases*}
      \circ_{j+k-1}\cdot(\circ_i\otimes_{j+k-1} \uu) & if $1 \leq i< j\leq n$, \\
       \circ_{j}\cdot(\uu\otimes_{j} \circ_{i-j+1}) & if $j \leq i<j+m$, \\
        \circ_j\cdot(\circ_{i-m+1}\otimes_{j} \uu) & if $j+m \leq i\leq n+m-1$.
    \end{cases*}
  \end{equation}
\end{itemize}
\end{definition}
In the associativity axiom we ignored the canonical isomorphisms of sources  of morphisms~(\ref{eq:associativity}), given by Lemma~\ref{lemma:otimes_associativity}.
The first and the third cases of~(\ref{eq:associativity}) are called the parallel associativity axioms and the second case is called the sequential associativity axiom.
\begin{example}
When $C$ is discrete, the operations  
$$X\otimes_i  Y = \displaystyle\bigoplus\limits_{c\in C} X^i_c \otimes Y^c$$
reduce to plain coproducts. The maps 
\[\begin{tikzcd}
{P_n \otimes_i P_m}  & P_{n+m-1}
	\arrow["\circ_i",from=1-1, to=1-2]
\end{tikzcd}\]
 are then composition maps of classical operads colored by the set $C$.
\end{example}
\begin{definition}
A unit for a $C$-operad $P$ in $\V$ is a collection of maps
\begin{equation}\label{eq:unit1}\begin{tikzcd}
I
& \colorop P(a;b),
	\arrow["u_f",from=1-1, to=1-2]
\end{tikzcd}\end{equation}
from the monoidal unit $I$ of $\V$, for any $f\colon a\rightarrow b$ of $C$, such that 
\begin{equation}\label{eq:unit2}\begin{tikzcd}
I
& \colorop P(b;c)\\
& \colorop P(a;d)
	\arrow["u_g",from=1-1, to=1-2]
 \arrow["\colorop P(f;h)",from=1-2, to=2-2]
 \arrow["u_{hgf}"',from=1-1, to=2-2]
\end{tikzcd}\end{equation}
commutes for all maps $f,g$ and $h$, for which the diagram makes sense, and the unit axioms are satisfied:
\begin{equation}\label{eq:unit_axioms}
\circ_1\cdot(u_f\otimes \uu) = \colorop P(\uu\,\cdots\,\uu;f) \text{ and } 
\circ_i\cdot \,(1\otimes u_f) = \colorop P(\uu\,\cdots\,f\,\cdots \,\uu;\uu).\end{equation}
A $C$-operad with a unit is called a unital $C$-operad.
\end{definition}

In equation~(\ref{eq:unit_axioms}) we ignored the unit isomorphisms $$X\otimes I\cong X\cong I\otimes X$$ of the monoidal category $\V$. The unit maps $u_f$ interpret morphisms of $C$ as actual unary operations of~$P$.

\begin{definition}
    A morphism of $C$-operads (resp.~unital $C$-operads) is a  morphism of the underlying $C$-collections, which preserves the compositions (resp.~compositions and unit). The category of $C$-operads in $\V$ is denoted by $\it{Op}^C(\V)$ and the category of unital $C$-operads in $\V$ is denoted by $\it{Op}_u^C(\V)$. 
\end{definition}
We can view the maps (\ref{eq:circ_i}) as a cowedge 
\begin{equation}\label{eq:circ_i_cowedge}
       \begin{tikzcd}
\big\{(P_n)^i_{c}\otimes \prescript{c}{}{(P_m)} & P_{m+n-1}\big\}_{c\in C}.
\arrow["\circ_i",from=1-1, to=1-2]
\end{tikzcd}\end{equation}
In the presence of a unit we have:
\begin{lemma}\label{lemma:cowedge_redundant}
    For a unital $C$-operad $P$, the condition on the maps (\ref{eq:circ_i_cowedge}) to form a cowedge is redundant.
\end{lemma}
\pf
Consider the diagram
\[\begin{tikzcd}
	{\colorop P(\Gamma\,c\,\Delta;x) \otimes\colorop P(\Lambda;d)}  \\
	{\colorop P(\Gamma\,c\,\Delta;x)\otimes I \otimes\colorop P(\Lambda;d)} \\
	{\colorop P(\Gamma\,c\,\Delta;x)\otimes\colorop P(d;c)\otimes\colorop P(\Lambda;d)} & {\colorop P(\Gamma\,c\,\Delta;x)\otimes\colorop P(\Lambda;c)} \\
	{\colorop P(\Gamma\,d\,\Delta;x)\otimes\colorop P(\Lambda;d)} & {\colorop P(\Gamma\Lambda\Delta;x)}.
	\arrow["\cong", from=1-1, to=2-1]
	\arrow["{\uu\otimes u_f \otimes \uu}", from=2-1, to=3-1]
	\arrow["{\circ_i\otimes \uu}", from=3-1, to=4-1]
	\arrow["{\uu\otimes \circ_1}", from=3-1, to=3-2]
	\arrow["{\circ_i}", from=4-1, to=4-2]
	\arrow["{\circ_i}"', from=3-2, to=4-2]
	\arrow["{\uu\otimes \colorop P(\uu\,\cdots\,\uu;f)}",curve={height=-24pt}, from=1-1, to=3-2]
	\arrow["{\colorop P(\uu\,\cdots\,f\,\cdots \,\uu;\uu)\otimes \uu}"',shift right=10, curve={height=50pt}, from=1-1, to=4-1]
\end{tikzcd}\]
The capital Greek letters stand for lists of objects of $C$. It is enough to show that the outer diagram commutes for every map $f\colon d\rightarrow c$ of $C$. But the top right part and left part commute by unitality and the bottom right part commutes by the associativity of~$\circ_i$.
\epf
The next proposition gives an alternative description of $C$-operads, using more general composition maps $\circ_i^f$. These maps combine two abstract operations along a morphism $f$ of~$C$, which connects their input and output.

\begin{proposition}\label{proposition:alternative_C_operad}  
A $C$-operad in $\V$ is equivalently a $C$-collection $P$ in~$\V$,
 together with maps $\circ_i^f$ in $\V$
\begin{equation*}\label{eq:circ_i^f}
\begin{tikzcd}
\colorop P (c_1 \,\cdots\, c_i \,\cdots \,c_n;c) \otimes \colorop P (d_1 \,\cdots\,d_m;d) & \colorop P (c_1 \,\cdots\, d_1 \,\cdots\,d_m \,\cdots \,c_n;c),
	\arrow["\circ^f_i",from=1-1, to=1-2]
\end{tikzcd}\end{equation*}
for any $n\geq1,m\geq0,1\leq i\leq n$, objects $c,c_1,\ldots, c_n,d,d_1,\ldots,d_m$ and $f\colon d\rightarrow c_i$ in~$C$,
which are
\begin{itemize}
    \item $C$-equivariant:
    for any two maps 
	\[\begin{tikzcd}
\colorop (\sigma\, f_1 \,\cdots \,f_n;\hspace{8pt}f)\colon \colorop (d_1  \,\cdots \,d_n;d)
& \colorop (c_1 \,\cdots\, c_n;c),
	\arrow[from=1-1, to=1-2]
\end{tikzcd}\]
\[\begin{tikzcd}
\colorop (\tau\, g_1 \,\cdots \,g_m;\hspace{8pt}g)\colon \colorop (b_1  \,\cdots \,b_m;b)
& \colorop (a_1 \,\cdots\, a_m;a),
	\arrow[from=1-1, to=1-2]
\end{tikzcd}\]
in $(\bS C)^{\it{op}}\times C$ and a map $h\colon a\rightarrow c_i,$

\begin{equation}\label{equation:C-equivariance}
\circ_i^h\cdot \big(\colorop P(\sigma\, f_1 \,\cdots \,f_n;\hspace{8pt}f)\otimes \colorop P(\tau\, g_1 \,\cdots \,g_m;\hspace{8pt}g)\big)=\colorop P ((\sigma\circ_i\tau)\, f_1 \,\cdots\, g_1 \,\cdots\,g_m \,\cdots \,f_n;\hspace{35pt}f)\cdot \circ_{\sigma(i)}^{f_ihg},
 \end{equation}
\item $C$-associative: for any $1\leq j\leq n, 0\leq m$ and $0\leq k$, 
\begin{equation}\label{equation:C-associativity}
    \circ^g_i\cdot (\circ^f_j\otimes \uu)=
    \begin{cases*}
      \circ^f_{j+k-1}\cdot(\circ^g_i\otimes \uu) & if $1 \leq i< j\leq n$, \\
       \circ^f_{j}\cdot(\uu\otimes \circ^g_{i-j+1}) & if $j \leq i<j+m$, \\
        \circ^f_j\cdot(\circ^g_{i-m+1}\otimes \uu) & if $j+m \leq i\leq n+m-1$. 
    \end{cases*}
  \end{equation}
\end{itemize}
\end{proposition}
\pf
    The proof is a matter of rewriting the data and axioms of a $C$-operad.
\epf
We now relate unital $C$-operads to symmetric substitudes of \cite{Day:substitution}. The correspondence is analogous to the classical correspondence of May operads and unital Markl operads, see \cite[Proposition 13]{markl2008operads}.

Substitudes were originally defined for a small $\V$-enriched category $C$. We consider only $\V$-categories which arise from ordinary small categories. That is,
if $\V$ has coproducts and the product $\otimes$ is distributive, any small category $C$ is $\V$-enriched via 
\begin{equation}\label{eq:enrichment}
    C(c,d):=\displaystyle\bigoplus_{c \xrightarrow{f} d}I,
\end{equation}
the coproduct of the monoidal unit $I$ of $\V$. 
The $\V$-functors $C\rightarrow \V$ and their $\V$-natural transformations then agree with ordinary functors $C\rightarrow \V$ and their natural transformations. We reformulate the definition of a symmetric $\V$-substitude:
\begin{definition}
Let $C$ be a small category, viewed as a $\V$-category. A $\V$-substitude $P$ is a functor 
\[\begin{tikzcd}
(\bS C)^{\it{op}} \times C 
  &  \V,
	\arrow["P",from=1-1, to=1-2]
\end{tikzcd}\]
 together with:
\begin{itemize}
\item a family of morphisms
\[\begin{tikzcd}
\colorop P(x_1 \,\cdots\,
x_n;x)\otimes \colorop P(x_{11} \,\cdots\,
x_{1m_1};x_1)\otimes \cdots\otimes \colorop P(x_{n1} \,\cdots\,
x_{nm_n};x_n)
  &  \colorop P(x_{11} \,\cdots\,
x_{nm_n};x)
	\arrow["\mu",from=1-1, to=1-2]
\end{tikzcd}\]
natural in $x_{11} ,\ldots,
x_{nm_n}$ and $x$, called substitution, and
\item a family of morphisms
   \[\begin{tikzcd}
C(x,y) &  \colorop P(x;y)
	\arrow["\eta",from=1-1, to=1-2]
\end{tikzcd}\]
natural in $x$ and $y$, called the unit,
\end{itemize}
subject to obvious equivariance, associativity and unit axioms (cf.~\cite{Day:substitution}).
\end{definition}
 Note, that naturality of the substitution in the connecting variables $x_1,\ldots,x_n$ is not assumed explicitly, but follows from the unit and the associativity axioms, similarly as in Lemma~\ref{lemma:cowedge_redundant}.
\begin{proposition}\label{prop:C-operads_and_Substitudes}
Let $C$ be a small category.
    The categories of $\V$-substitudes over $C$ and unital $C$-operads in $\V$ are isomorphic.
\end{proposition}
\pf
The underlying collections of the two structures are the same. Since $C(c,d)$ is defined by (\ref{eq:enrichment}), the data of a unit for both structures also agree. Due to Lemma~\ref{lemma:cowedge_redundant} we are left to show the correspondence of the compositions $\circ_i$ and the substitution $\mu$, which is standard. The operation $\mu$ is built up from consecutive $\circ_i$ operations and conversely, the operations $\circ_i$ are obtained from $\mu$ by inserting units into the inputs of~$\mu$. 
\epf

\begin{remark}\label{remark:decomposition}
As noted in \cite[Remark 5.2.1]{Batanin_White:substitudes}, $C$-collections are the same as $C\bS$-modules of~\cite{petersen2013operad}, up to variance and order of variables, and thus also the same as \hbox{$\mathbb{V}$-sequences} of~\cite{Ward:Massey}. A \hbox{$\V$-substitude} over a small category $C$ corresponds to a $C$-colored operad in $\V$ of \cite{petersen2013operad}, and hence to a $C$-operad of \cite{Ward:Massey}, when $C$ is a groupoid. 
However, it seems that \cite{Ward:Massey} works only with skeletal groupoids, i.e.~groupoids with $C(a,b)=\emptyset$ when $a\neq b$. 
For us, it is essential to allow morphisms between different objects of the coloring groupoid, since we encode Markl \hbox{$\bO$-operads}, which are $\bO_{\it{iso}}$-collections in the first place. An example of an operadic category with isomorphisms between different objects is the operadic category $\texttt{Tr}$ of rooted trees (cf.~Section~3 and Example~4.8 of~\cite{Batanin_Markl:kodu2022}, and the introduction of \cite{Batanin_Markl_Obrad:models}). Algebras over the terminal $\texttt{Tr}$-operad are classical symmetric operads.
\end{remark}

\section{Categorical operads and internal operads}\label{section:internal_operads}
\begin{figure}[H]
    \centering
    \includegraphics{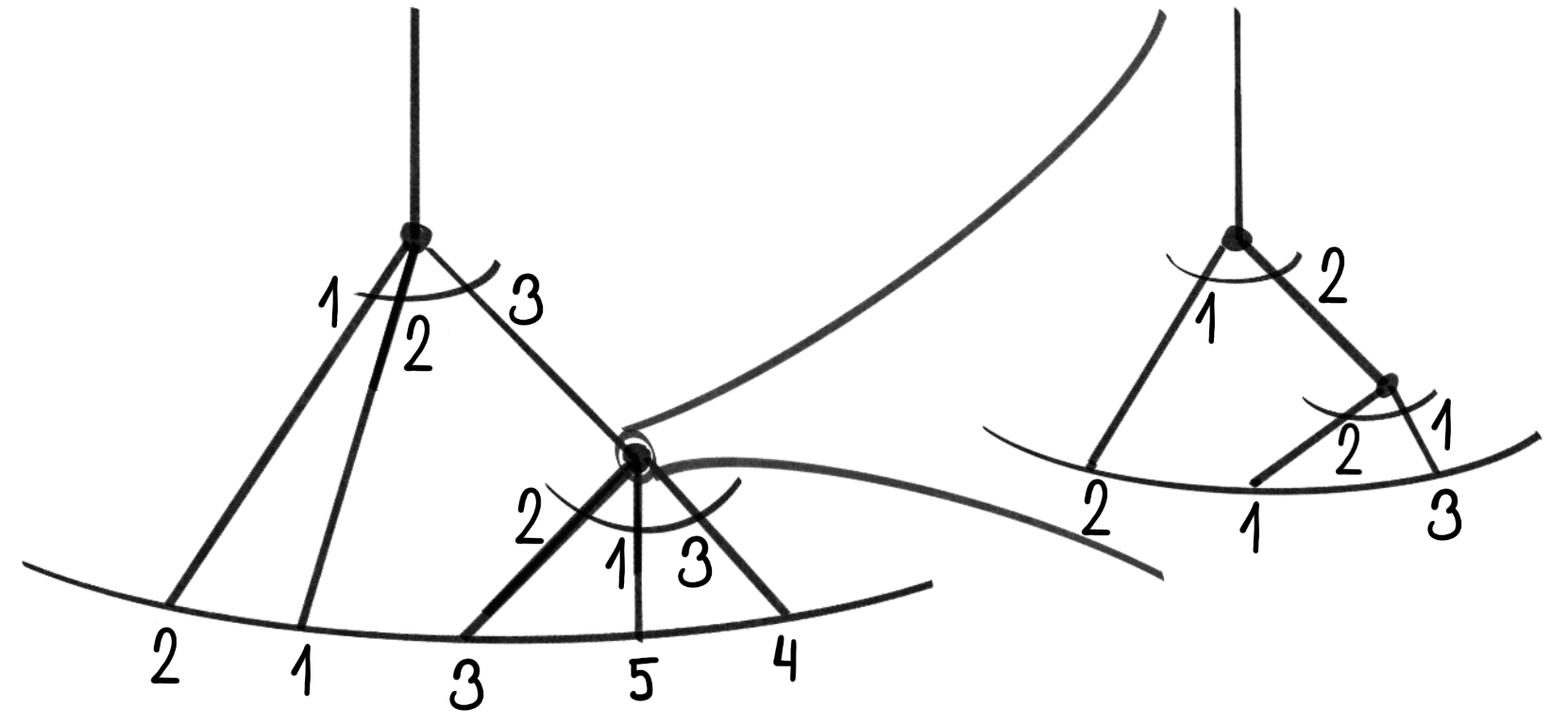}
    \caption{Ordered trees and substitution.}
    \label{fig:tree_monad}
\end{figure}

Recall from \cite[Section 9]{Batanin_Berger:polynomial} the polynomial monad $T$ on the slice category $\textit{Set}/\mathbb{N}$ (the category of families of sets indexed by natural numbers), whose category of algebras is isomorphic to the category of symmetric operads in $\textit{Set}$. The monad is generated by the polynomial 
 $$\mathbb{N} \xleftarrow{s} \texttt{ORTrees}^* \xrightarrow{p} \texttt{ORTrees} \xrightarrow{t} \mathbb{N},$$
 where 
\texttt{ORTrees} is a set of isomorphism classes of rooted trees, which have a linear order
on the set of leaves and also on the sets of incoming edges for each vertex of the tree.
$\texttt{ORTrees}^*$ stands for ordered rooted trees with one marked vertex. The map $p$ forgets the marking, the target map $t$ sends a tree to the number of its leaves, and the source map $s$ gives a number of inputs of the marked vertex. The monad $T$ is defined by the formula
$$ TX_n=\displaystyle\sum_{a\in t^{-1}(n)}\displaystyle\prod_{b\in p^{-1}(a)}X_{s(b)}.$$
The multiplication is given by substituting a tree into the marked vertex,  as in Figure~\ref{fig:tree_monad}.
The polynomial monad $T$ generates a cartesian 2-monad (also denoted by $T$) on the \hbox{2-category} $\textit{Cat}/\mathbb{N}$ of small categories over the discrete category of natural numbers.
\begin{definition}
A unital categorical operad is a pseudo-$T$-algebra. The category of unital categorical operads is denoted by $\it{Op}_u$.
\end{definition}
\begin{definition}\label{definition:internal_operad}
Let $\OO$ be a unital categorical operad. An internal unital operad $P$ in~$\OO$ is a lax morphism from the terminal pseudo-$T$-algebra to $\OO$. The category of internal unital operads in~$\OO$ and $T$-natural transformations will be denoted by $\it{IOp}_u(\OO)$.
\end{definition}
We consider variants of the monad $T$ which describe non-unital and non-symmetric categorical operads.
Let $T_{nu}$ be the polynomial monad generated by the polynomial 
$$\mathbb{N} \xleftarrow{s} \texttt{ORTrees}_+^* \xrightarrow{p} \texttt{ORTrees}_+ \xrightarrow{t} \mathbb{N},$$
 where 
$\texttt{ORTrees}_+$ is the set of isomorphism classes of ordered rooted trees, without the tree~$|$ (the free living edge). Consequently, algebras of the monad $T_{nu}$ have no units for their composition operations. 
Let $T_{\it{ns}}$ be generated by the polynomial 
$$\mathbb{N} \xleftarrow{s} \texttt{PRTrees}_+^* \xrightarrow{p} \texttt{PRTrees}_+ \xrightarrow{t} \mathbb{N},$$
 where 
$\texttt{PRTrees}_+$ is the set of planar rooted trees, see \cite[Section 14]{Batanin_Berger:polynomial}, without the tree~$\,|\,$.
\begin{definition}
A categorical operad is a pseudo-$T_{nu}$-algebra. The category of categorical operads is denoted by $\it{Op}$. A non-symmetric categorical operad is a pseudo-$T_{\it{ns}}$-algebra. The category of non-symmetric categorical operads is denoted by $\it{Op}_{\it{ns}}$.
\end{definition}
\begin{definition}
An internal operad in a categorical operad $\OO$ is a~lax morphism from the terminal categorical operad into $\OO$. The category of internal operads in $\OO$ will be denoted~$\it{IOp}(\OO)$. An internal non-symmetric operad in a non-symmetric categorical operad $\OO$ is a lax morphism from the terminal non-symmetric categorical operad into $\OO$. The~category of internal non-symmetric operads in $\OO$ will be denoted $\it{IOp}_{\it{ns}}(\OO)$.
\end{definition}

The identity functor $\uu_{\it{Cat}/\mathbb{N}}$ on the category $\it{Cat}/\mathbb{N}$ is a monad induced by a constant polynomial and $\uu_{\it{Cat}/\mathbb{N}}$-algebras are just objects of $\it{Cat}/\mathbb{N}$. An internal algebra $x$ in a $\uu_{\it{Cat}/\mathbb{N}}$-algebra $X$ is a collection of object $\{x_n \in X(n)\}_{n\geq0}$. We will refer to it as to internal collection in $X$. The category of internal collections in $X$ will be denoted by $\it{IColl}(X)$.

Any planar tree can be regarded as an ordered tree, for which the global order on its leaves agrees with the order induced by local orders of vertices and planarity.  We have inclusions 
\[\begin{tikzcd}
	\NN & {\texttt{PRTrees}_+} &{\texttt{ORTrees}_+}& {\texttt{ORTrees}}. 
	\arrow[hook, from=1-1, to=1-2]
	\arrow[hook, from=1-2, to=1-3]
	\arrow[hook, from=1-3, to=1-4]
\end{tikzcd}\]
The leftmost inclusion interprets a natural number $n$ as an $n$-corolla (i.e.~a planar rooted tree with one vertex and $n$ leaves). The inclusions induce maps of polynomial monads and hence functors between pseudo-algebras:
\[\begin{tikzcd}
	\Cat/\NN & \it{Op}_{\it{ns}} &\it{Op}& \it{Op}_u. 
	\arrow[from=1-2, to=1-1]
	\arrow[from=1-3, to=1-2]
	\arrow[from=1-4, to=1-3]
\end{tikzcd}\]
In other words, we can view any unital categorical operad as a non-unital/non-symmetric categorical operad. For any unital categorical operad $\OO$ we have the forgetful functors
\[\begin{tikzcd}
	\it{IColl}(\OO) & \it{IOp}_{\it{ns}}(\OO) &\it{IOp}(\OO)& \it{IOp}_u(\OO). 
	\arrow["U_1"',from=1-2, to=1-1]
	\arrow["U_2"',from=1-3, to=1-2]
	\arrow["U_3"',from=1-4, to=1-3]
\end{tikzcd}\]
In Section \ref{section:free} we study the left adjoint to the composite $U=U_1\circ U_2\circ U_3.$

A categorical operad can be presented as a collection of categories $\OO(n)$, for each $n\geq0$, together with 
\begin{itemize}
    \item  functors
\[\begin{tikzcd}
	\OO(n) & \OO(n),
	\arrow["{\sigma(-)}", from=1-1, to=1-2]
\end{tikzcd}\]
    for any $n\geq1$ and permutation $\sigma \in \Sigma_n$, making each $\OO(n)$ a non-strict left $\Sigma_n$-module, i.e.~there are coherent natural isomorphism $\tau(\sigma(-) ) \cong \tau\sigma(-) $ and $\uu_n(-)\cong \uu_{\OO(n)}$,
    \item partial compositions functors
\[\begin{tikzcd}
\OO(n) \times \OO(m) &  \OO(m+n-1)
	\arrow["\odot_i", from=1-1, to=1-2]
\end{tikzcd}\]
for $n>0,m\geq0, $ and $ i \in \{1,\ldots,n\}$,
\end{itemize}
satisfying classical associativity and equivariance axioms (see \cite[Definition 11]{markl2008operads}), only up to coherent isomorphisms.

Consequently, an internal operad $A$ in $\OO$ is a collection of objects $A_n \in \OO(n)$, for each $n\geq0$, with the following additional structure.
\begin{itemize}
   \item Each $A_n$ is equipped with an internal $\Sigma_n$-action: For every $\sigma\in \Sigma_n$ there is a map
    \[\begin{tikzcd}
	\sigma (A_n) &  A_n
	\arrow["{\overline{\sigma}}", from=1-1, to=1-2]
\end{tikzcd}\]
in the category $\OO(n)$, such that
    \[\begin{tikzcd}
	 \sigma(\tau (A_n)) & \sigma (A_n)\\
	  \sigma \tau(A_n)& A_n
	\arrow["\sigma(\overline{\tau})", from=1-1, to=1-2]
 \arrow["\overline{\sigma}", from=1-2, to=2-2]
	\arrow["\overline{\sigma\tau}", from=2-1, to=2-2]
	\arrow["{\cong}", from=1-1, to=2-1]
\end{tikzcd}\]
commutes and
\[\begin{tikzcd}
	\uu_n(A_n) &  A_n
	\arrow["{\overline{\uu_n}}", from=1-1, to=1-2]
 \end{tikzcd}\]
is the natural isomorphisms of the non-strict $\Sigma_n$-action on $\OO(n)$.
\item For every $n>0, m\geq 0$ and $1\leq i \leq n$, there are internal composition maps
\[\begin{tikzcd}
A_n\odot_i A_m  &  A_{m+n-1}
	\arrow["\bullet_i", from=1-1, to=1-2]
\end{tikzcd}\]
in the category $\OO(m+n-1)$.

\end{itemize}
The structure maps further satisfy internal versions of classical operad axioms:

\begin{itemize}
\item \textit{Associativity.}
 For any $1\leq j\leq n, 0\leq m$ and $0\leq k$, 
\begin{equation}\label{eq:internal_associativity}
    \bullet_i\cdot (\bullet_j\odot_i \uu)=
    \begin{cases*}
      \bullet_{j+k-1}\cdot(\bullet_i\odot_{j+k-1} \uu) & if $1 \leq i< j\leq n$, \\
       \bullet_{j}\cdot(\uu\odot_{j} \bullet_{i-j+1}) & if $j \leq i<j+m$, \\
        \bullet_j\cdot(\bullet_{i-m+1}\odot_{j} \uu) & if $j+m \leq i\leq n+m-1$. 
    \end{cases*}
  \end{equation}
In formulas (\ref{eq:internal_associativity}) we ignored the structure isomorphisms of the operad $\OO$. For example, the middle case thus means the commutativity of the diagram 
   \[\begin{tikzcd}
	{(A_{n} \odot_j  A_{m})\odot_i A_{k}} & {} & {A_{n+m-1}\odot_i A_{k}} \\
	{A_{n} \odot_j  (A_{m}\odot_{i-j+1} A_{k})} \\
	{A_{n}\odot_j A_{m+k-1}} & {} & {A_{n+m+k-2}}.
	\arrow["{\cong}", from=1-1, to=2-1]
	\arrow["{\uu\odot_j \bullet_{i-j+1}}", from=2-1, to=3-1]
	\arrow["{\bullet_j\odot_i  \uu}", from=1-1, to=1-3]
	\arrow["{\bullet_j}", from=3-1, to=3-3]
	\arrow["{\bullet_i}", from=1-3, to=3-3]
\end{tikzcd}\]
\item \textit{Equivariance.} For $\sigma\in\Sigma_n$, $\tau\in\Sigma_m$, $1\leq i\leq n$ and $\pi=\sigma\circ_i\tau$, the following diagram commutes.

\begin{equation}\label{eq:internal_equivariance_axioms}\begin{tikzcd}
	{\pi(A_{n} \odot_{\sigma(i)}  A_{m})} & {} & {\pi(A_{n+m-1})}\\
	{\sigma(A_{n}) \odot_{i}  \tau(A_{m})} && {} \\
	{A_{n} \odot_{i}  A_{m}} & {} & {A_{n+m-1}}
	\arrow["{\overline{\sigma}\odot_{i}  \overline{\tau}}", from=2-1, to=3-1]
	\arrow["{\pi(\bullet_{\sigma(i)})}", from=1-1, to=1-3]
	\arrow["{\cong}", from=1-1, to=2-1]
	\arrow["{\bullet_{i}}", from=3-1, to=3-3]
	\arrow["{\overline{\pi}}", from=1-3, to=3-3]
\end{tikzcd}\end{equation}
\end{itemize} The isomorphisms in the diagrams are the structure isomorphisms of the operad~$\OO$.
Non-symmetric internal operads admit a similar presentations without the structure maps and axioms involving internal $\Sigma$-actions.
In the unital case, a unital categorical operad~$\OO$ further has a unit $U\in \OO(1)$, satisfying the unit axioms up to coherent natural isomorphisms. A~unital internal operad then has an internal unit, i.e.~a morphism 
\[\begin{tikzcd}
U  &  A_{1}
	\arrow["u", from=1-1, to=1-2]
\end{tikzcd}\] in the category $\OO(1)$, 
such that for $1\leq i\leq n$, the following diagrams commute.
\begin{equation}\label{eq:internal_unit_axioms}
    \begin{tikzcd}
	{A_{n}  \odot_{i}  U} & {A_{n}  \odot_{i}  A_1} \\
	{} & {A_{n} }
	\arrow["{\cong}"', from=1-1, to=2-2]
	\arrow["{\uu  \odot_{i}  u}", from=1-1, to=1-2]
	\arrow["{\bullet_i}", from=1-2, to=2-2]
\end{tikzcd}
\hspace{1cm}
\begin{tikzcd}
	 {U  \odot_{1}  A_{m}} & {A_1\odot_{1}  A_{m}} \\
	  & {A_{m}}
	\arrow["{\bullet_{1}}", from=1-2, to=2-2]
	\arrow["{u \odot_{1}  \uu }", from=1-1, to=1-2]
	\arrow["{\cong}"', from=1-1, to=2-2]
\end{tikzcd}
\end{equation}
The following example is an analogue of Lemma 9.1 of \cite{Batanin:eckmann}.
\begin{example}\label{example:1}
Let $(\V,\otimes,I)$ be a symmetric monoidal category. Define $\OO(n):=\V$, with the trivial $\Sigma$-actions, $\odot_i:= \otimes$ and $U:=I$. The non-strict associativity of $\otimes$ makes $\OO$ a unital categorical operad. The category of internal operads in $\OO$ is isomorphic to the category of classical operads in $\V$.
\end{example}
\begin{remark} In the context of Example~\ref{example:1}, one can think of the structure operations~$\odot_i$ of~$\OO$ as a generalisation of the product $\otimes$ of $\V$. When composing three (or more) operations of a classical operad $P$, we loose the shape of the source $P_n\otimes P_m \otimes P_k$.
On the other hand, when defining operads in a categorical operad $\OO$, the source $(P_n\odot_j P_m) \odot_i P_k$ of multiple compositions retains its geometrical shape, which makes the compositions easier to handle, for example when describing free operads.
\end{remark}

\section{Categorical operad $\SSS$ for $C$-operads}\label{section:SSS}
In this section, we describe a categorical operad $\SSS$, for a small category $C$ and cocomplete symmetric monoidal closed category $\V$, and prove that internal operads in $\SSS$ are \hbox{$C$-operads}.
Recall the notation $\SSS(n)$, for the category of functors 
\[\begin{tikzcd}
\underbrace{C^{\it{op}} \times \cdots \times C^{\it{op}}}_{n\textit{-times}} \times C  &  \V,
	\arrow["X",from=1-1, to=1-2]
\end{tikzcd}\]
 and the operations 
 \[\begin{tikzcd}
	\SSS(n)\times \SSS(m) & \SSS(n+m-1)
 \arrow["\otimes_i", from=1-1, to=1-2]
\end{tikzcd}\]
of Definition \ref{definition:odot_i}.
Let us describe the left $\Sigma_n$-module structure on $\SSS(n)$.
Let  $\underline{n}$ denote the set $\{1,\ldots,n\}$, regarded as a discrete category.
There is a canonical right $\Sigma_n$-module structure on the functor category $\Cat(\underline{n},C)$ given by precomposition:
\[\begin{tikzcd}
	\underline{n} & C. \\
	\underline{n} & 
	\arrow["\sigma", from=2-1, to=1-1]
	\arrow["f", from=1-1, to=1-2]
	\arrow["{ f\circ\sigma }"', dashed, from=2-1, to=1-2]
\end{tikzcd}\]
We transfer this right $\Sigma_n$-module structure from $\Cat(\underline{n},C)$ to $C^{\times n}$ by
\[\begin{tikzcd}
	C^{\times n} & C^{\times n} \\
	{\Cat(\underline{n},C)} & {\Cat(\underline{n},C)}
	\arrow["\cong"', from=1-1, to=2-1]
	\arrow["{(-)\circ\sigma}", from=2-1, to=2-2]
	\arrow["\cong"', from=2-2, to=1-2]
	\arrow["{(-)\sigma}", dashed, from=1-1, to=1-2]
\end{tikzcd}\]
\[\begin{tikzcd}
	{\underline{c}=(c_1,\ldots,c_n)} & {(c_{\sigma(1)},\ldots,c_{\sigma(n)})=\underline{c}\sigma}
	\arrow["{(-) \sigma}", maps to, from=1-1, to=1-2]
\end{tikzcd}\]
using the obvious isomorphism 
$$C^{\times n}\cong \Cat(\underline{n},C).$$
The left $\Sigma_n$-module structure on $\SSS(n)$ is given by 
\[\begin{tikzcd}
	{(C^{\it{op}})^{\times n}\times C} & \V.\\
	{(C^{\it{op}})^{\times n}\times C} & {}
	\arrow["X", from=1-1, to=1-2]
	\arrow["{(-)\sigma \times \uu}", from=2-1, to=1-1]
\arrow["{\sigma(X)}"', dashed, from=2-1, to=1-2]
\end{tikzcd}\]
\begin{lemma}\label{lemma:otimes_equivariance}
The operations $\otimes_i$ are equivariant.
\end{lemma}
\pf
For any $X \in \SSS(n)$, $Y \in \SSS(m)$, $\sigma \in \Sigma_n$, $\tau \in \Sigma_m$, $f\colon c \rightarrow d \in C$, and $i\in \{1,\ldots,n\}$, there is a strict equality
\begin{equation*}\label{equality:equivariance}
\sigma (X)^i_c \otimes \prescript{d}{}{\tau (Y)} = (\sigma\circ_i \tau)(X^{\sigma(i)}_c \otimes  \prescript{d}{}{Y}).
\end{equation*}
  Further, the map
  \[\begin{tikzcd}
	\SSS(n+m-1) & & \SSS(n+m-1)
 \arrow["(\sigma\circ_i\tau)(-)", from=1-1, to=1-3]
\end{tikzcd}\]
 is an isomorphism, so it commutes with colimits, and thus $(\sigma\circ_i\tau)(X\otimes_{\sigma(i)} Y)$ is the colimit of the diagram 
\[\begin{tikzcd}
	{\displaystyle\bigoplus_{f\colon d \rightarrow c}(\sigma\circ_i\tau)(X^{\sigma(i)}_c \otimes \prescript{d}{}{Y})} &&& {\displaystyle\bigoplus_{b\in C} (\sigma\circ_i\tau)(X^i_b \otimes \prescript{b}{}{Y}).}
	\arrow["{\displaystyle\bigoplus_{f\colon d \rightarrow c}(\sigma\circ_i\tau)(X^{\sigma(i)}_f\otimes \uu)}", curve={height=-30pt}, from=1-1, to=1-4]
	\arrow["{\displaystyle\bigoplus_{f\colon d \rightarrow c}(\sigma\circ_i\tau)(\uu \otimes \prescript{f}{}{Y})}"', curve={height=30pt}, from=1-1, to=1-4]
\end{tikzcd}\]
One obtains the canonical natural equivariance isomorphism
\[\sigma (X)\otimes_{i} \tau (Y) \cong (\sigma\circ_i\tau)(X\otimes_{\sigma(i)} Y).
\]
\epf

The associativity of $\otimes_i$ was shown in Lemma \ref{lemma:associativity_of_X^i_otimes_Y}. Since the isomorphisms constructed in Lemma \ref{lemma:otimes_associativity} and Lemma \ref{lemma:otimes_equivariance} are canonical, we have proved: 
\begin{proposition}\label{proposition:C_V_is_an_operad}
The categorical left $\Sigma_n$-modules $\SSS(n)$ together with operations $\otimes_i$ form a categorical operad.
\end{proposition}
The operad $\SSS$ is in fact unital. The unit for $\otimes_i$ is the functor $C(-,-) \in \SSS(1)$, defined by~(\ref{eq:enrichment}).
An internal operad $A$ in the categorical operad $\SSS$ consists of a collection \hbox{$A_n\in \SSS(n)$} for each $n\geq 0$, i.e.~a collection of functors 
\[\begin{tikzcd}
\underbrace{C^{\it{op}} \times \cdots \times C^{\it{op}}}_{n\textit{-times}} \times C  & \V,
	\arrow["A_n",from=1-1, to=1-2]
\end{tikzcd}\]
together with natural transformations 
\[\begin{tikzcd}
\sigma (A_n)   & A_n.
	\arrow["\overline{\sigma}",from=1-1, to=1-2]
\end{tikzcd}\]
By definition,
$$\sigma (A_n)\colorop(x_1 \,\cdots\,
x_n;x)=\colorop A_n(x_{\sigma(1)} \,\cdots\,
x_{\sigma(n)};x),$$
and so the components of $\overline{\sigma}$ are 
\[\begin{tikzcd}
\colorop A_n (x_{\sigma(1)} \,\cdots\,
x_{\sigma(n)};x)
&&&
\colorop A_n (x_1 \,\cdots\,x_n;x).
	\arrow["\overline{\sigma}\colorop (x_1 \,\cdots\,
x_n;x)",from=1-1, to=1-4]
\end{tikzcd}\]
Further, there are natural transformations 
\[\begin{tikzcd}
A_n \otimes_i A_m & A_{n+m-1},
	\arrow["\bullet_i",from=1-1, to=1-2]
\end{tikzcd}\]
the unit is a natural transformation
\[\begin{tikzcd}
C(x,y) & \colorop A_1(x;y),
	\arrow["u_{xy}",from=1-1, to=1-2]
\end{tikzcd}\]
and they satisfy the internal operad axioms (\ref{eq:internal_associativity}), (\ref{eq:internal_equivariance_axioms}), and~(\ref{eq:internal_unit_axioms}).

\begin{theoremA}
\label{proposition:iso_of_int_operads_and_C_operads}
The category of internal operads (resp.~unital internal operads) in the categorical operad $\SSS$ (resp.~unital categorical operad $\SSS$), constructed above, is isomorphic to the category of $C$-operads (resp.~unital $C$-operads) in $\V$:
$$
\it{IOp}(\SSS)\cong \it{Op}^C(\V) 
\text{ and }
\it{IOp}_u(\SSS)\cong \it{Op}_u^C(\V). 
$$
\end{theoremA}
\pf
 The data of the underlying collection of a $C$-operad and of an internal operad agree. The internal actions of an internal operad 
\[\begin{tikzcd}
 \sigma\colorop A_n(x_1 \,\cdots\,
x_n;x)
& 
A_n\colorop (x_1 \,\cdots\,
x_n;x)
	\arrow["\overline{\sigma}",from=1-1, to=1-2]
\end{tikzcd}\]
translate as right $\Sigma_n$-actions of a $C$-operad
\[\begin{tikzcd}
 P_n\colorop (x_{\sigma(1)} \,\cdots\,
x_{\sigma(n)};x)
&
P_n \colorop (x_1 \,\cdots\,x_n;x).
	\arrow[" (-)\sigma^{-1}",from=1-1, to=1-2]
\end{tikzcd}\]
By definition, the operations $\circ_i$ of a $C$-operad and $\bullet_i$ of an internal operad are the same. The units for both structures agree and
the equivariance and associativity axioms can be compared directly.
\epf
When $C$ is a small $\V$-enriched category, the categories $\hat{C_{\V}}(n)$ of $\V$-functors and their $\V$-natural transformations again form a categorical operad. The composition operations $\otimes_i$ are defined using $\V$-enriched colimits. Theorem A suggest the following definition of $C$-operads for a $\V$-enriched category $C$.
\begin{definition}\label{definition:nu_C_operad_V}
    Let $C$ be a small $\V$-enriched category. A $C$-operad $P$ is an internal operad in the categorical operad $\hat{C_{\V}}$.
\end{definition}
An example of a $\Vect$-enriched $C$-operad is given in Example \ref{example:dgas}.
\section{Free internal operads}\label{section:free}
We recall some combinatorial material and terminology of \cite{Batanin_Markl:kodu2022}, which we will use to describe free internal operads.

\subsection*{The operadic category $\Delta$}

Let $\Delta$ be the category of finite ordered sets $ \underline{n}  := \{1, \ldots, n\}$ and order preserving maps. It is an operadic category in the sense of \cite[Section~1]{Batanin_Markl:kodu2022}. Recall that a map $\phi\colon \underline{m} \rightarrow
\underline{n}$ of $\Delta$ is \textit{elementary} if it has precisely one
non-trivial fiber, i.e.~there is a unique index $i \in \underline{n}$, such that $|\phi^{-1}(i)|\neq \underline{1}$. We write $(\phi,i)$ for an
elementary morphism $\phi$ with a non-trivial fiber over $i$. 
Two composable elementary morphisms $$\underline{n} \xrightarrow{(\psi,j)} \underline{m}\xrightarrow{(\phi,i)}\underline{k}$$ have \textit{disjoint
fibers} if $\phi(j)\neq i$.

\begin{figure}[H]
    \centering
    \includegraphics{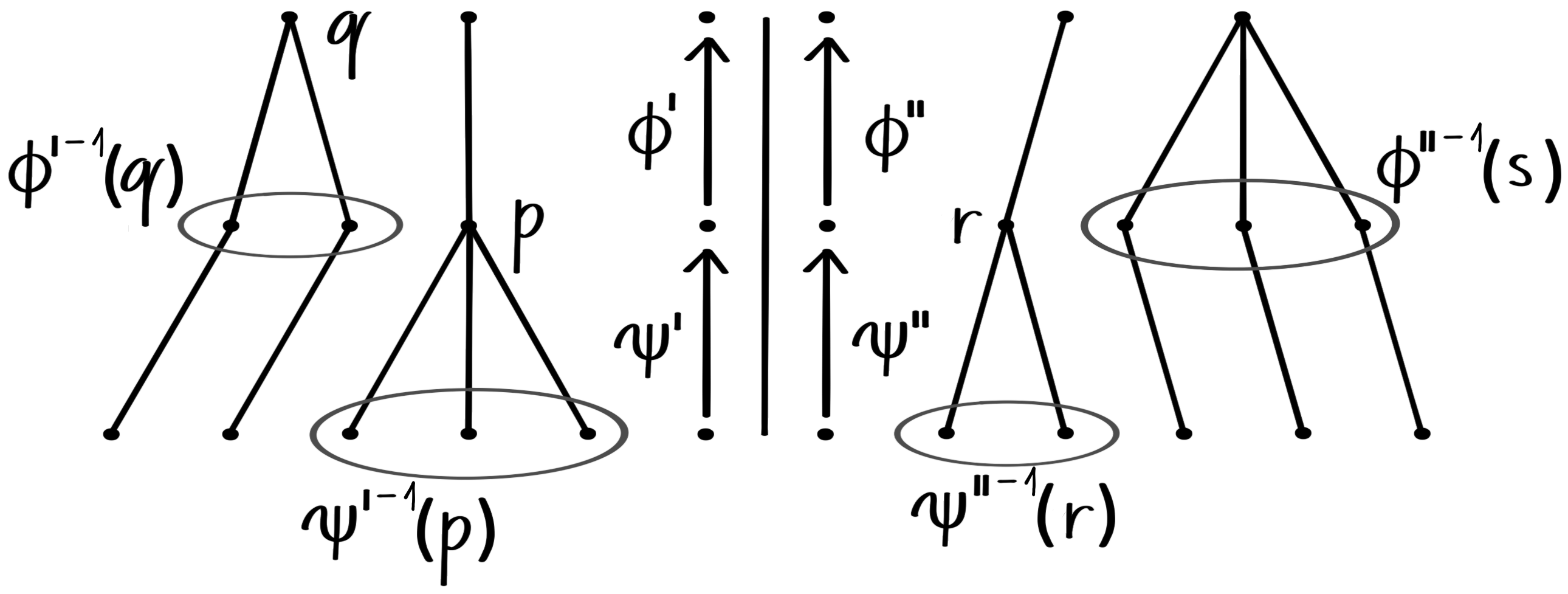}
    \caption{Pair of morphisms with disjoint fibers.}
    \label{fig:fibers}
\end{figure}

Consider the following commutative diagram of elementary morphisms with
disjoint fibers.
\begin{equation}\label{diagram:square}
    \begin{tikzcd}
	& {\underline{m}'} \\
	{\underline{n}} && {\underline{k}} .\\
	& {\underline{m}''}
	\arrow["{(\psi',p)}", from=2-1, to=1-2]
	\arrow["{(\phi',q)}", from=1-2, to=2-3]
	\arrow["{(\psi'',r)}"', from=2-1, to=3-2]
	\arrow["{(\phi'',s)}"', from=3-2, to=2-3]
\end{tikzcd}
\end{equation}
In this situation the non-trivial fibers of $\psi'$ and $\phi''$ are the same (cf.~Figure~\ref{fig:fibers}): $$\psi'^{-1}(p)=\phi''^{-1}(s),$$ and similarly $$\psi''^{-1}(r)=\phi'^{-1}(q).$$

\subsection*{Free construction}

Fix a unital categorical operad $\OO$, for which every $\OO(n)$ is cocomplete and the operations $\odot_i$ commute with colimits in both variables.
Recall the restriction functors 
\begin{equation}\label{equation:U}
\begin{tikzcd}
	\it{IColl}(\OO) & \it{IOp}_{\it{ns}}(\OO) &\it{IOp}(\OO)& \it{IOp}_u(\OO). 
	\arrow["U_1"',from=1-2, to=1-1]
	\arrow["U_2"',from=1-3, to=1-2]
	\arrow["U_3"',from=1-4, to=1-3]
\end{tikzcd}
\end{equation}
of Section~\ref{section:internal_operads} and denote their composite by $U$. A left adjoint $F \dashv U$ exists and can be computed by the methods of~\cite{Batanin_Berger:polynomial}, using the relative internal algebra classifiers. We~describe~$F$ by constructing the left adjoints $F_i \dashv U_i$ for each functor of the sequence~(\ref{equation:U}).

\begin{figure}[H]
    \centering
    \includegraphics{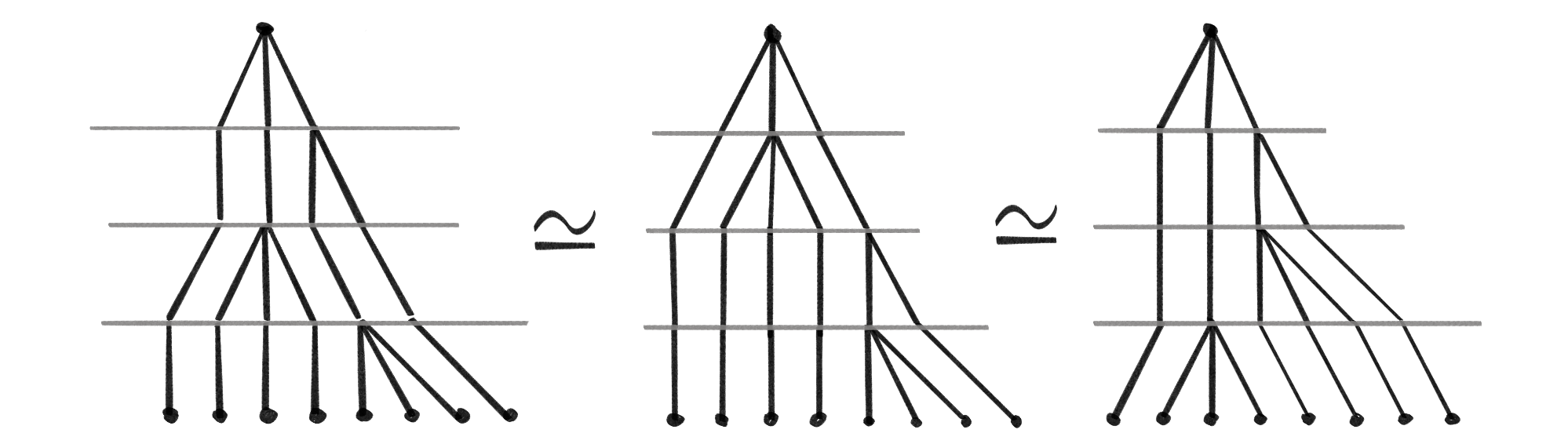}
    \caption{The groupoid of leveled trees.}
    \label{fig:leveled_trees}
\end{figure}

For the construction of $F_1 \dashv U_1$ we use a modification of the construction of \cite{Batanin_Markl:kodu2021} for the operadic category $\Delta$, which produces the free operad as a colimit of a functor from the groupoid of leveled planar rooted trees. We replace the monoidal products $\otimes$ in the construction with the operations $\odot_i$.
\begin{definition}
	A \textit{leveled planar rooted tree} $T$ is a sequence of composable elementary morphisms $$\underline{n}_1
	\xrightarrow{(\tau_1,i_1)} \cdots
	\xrightarrow{(\tau_{k-1},i_{k-1})}
	\underline{n}_k \xrightarrow{!} \underline{1}$$ in $\Delta$. Its \textit{fiber sequence} $f_1 ,\ldots, f_k$ is a sequence of non-trivial fibers of the
	morphisms. The height of a tree is the number of morphisms in the sequence.
\end{definition}
Let us define a groupoid $\lT$ of leveled planar rooted trees. The objects are leveled planar rooted trees and two trees
\begin{equation}\label{isomorphic-trees}
  \begin{split}
    \T' = \underline{n_1}
	\xrightarrow{} \cdots
	\xrightarrow{}
	\underline{n} \xrightarrow{}
	\underline{m}'\xrightarrow{} \underline{k}\xrightarrow{}\cdots \xrightarrow{} \underline{1} \\
 \T''= \underline{n_1}
	\xrightarrow{} \cdots
	\xrightarrow{}
	\underline{n} \xrightarrow{}
	\underline{m}''\xrightarrow{} \underline{k}\xrightarrow{}\cdots \xrightarrow{} \underline{1}
  \end{split}
\end{equation}
 are isomorphic if their sequences differ by a commutative
square of elementary morphisms with disjoint fibers as in (\ref{diagram:square}). Three isomorphic leveled planar rooted trees are depicted in Figure~\ref{fig:leveled_trees}.
Let $\lT(n)$ denote the full subcategory of trees with $\underline{n}_1=\underline{n}$, so the collection of $\lT(n)$ is an object in $\it{Cat}/\mathbb{N}$.
The maps 
\[\begin{tikzcd}
\lT(n) \times \lT(m)  & \lT(m+n-1)
	\arrow["\odot_i",from=1-1, to=1-2]
\end{tikzcd}\]
are given by grafting of leveled trees, producing a (strict) non-unital non-symmetric categorical operad $\lT$.
Let $$\T =\underline{n}_1\xrightarrow{(\tau_1,i_1)} \cdots \xrightarrow{(\tau_k,i_k)} \underline{1}$$ be a tree with fiber sequence
$t_1,\ldots, t_k$. For an internal collection $X$ in $\OO$ we define a functor
\[\begin{tikzcd}
\lT  & \OO
	\arrow["\hat{X}",from=1-1, to=1-2]
\end{tikzcd}\]
by the formula
$$\hat{X}(\T)= (\cdots\, (X_{t_k}\odot_{i_{k-1}} X_{t_{k-1}})\odot_{i_{k-2}}\cdots )\odot_{i_{1}}X_{t_1}.$$
The level-interchange isomorphisms of $\lT$ are mapped by $\hat{X}$ to the parallel associativity isomorphisms of $\OO$. Further, we have
$$\hat{X}(\T\odot_i \sS)\cong \hat{X}(\T)\odot_i \hat{X}(\sS).$$
The isomorphism is given by the repeated use of the sequential associativity isomorphisms of $\OO$.

\begin{definition}
 Let $X$ be an internal collection in $\OO$. We define the functor $F_1$ by 
 \[\begin{tikzcd}
F_1(X):= \it{colim} \big(\lT & \OO\big).
	\arrow["\hat{X}",from=1-1, to=1-2]
\end{tikzcd}\]
\end{definition}
By assumption the operations $\odot_i$ of the operad $\OO$ commute with colimits, and thus
\begin{equation}\label{eq:composition_in_F_1}
F_1(X)_n\odot_i F_1(X)_m = (\underset{\lT(n)}{\it{colim}}\hspace{5pt} \hat{X}_n)\odot_i(\underset{\lT(m)}{\it{colim}}\hspace{5pt} \hat{X}_m)\cong \underset{\lT(n)\times \lT(m)}{\it{colim}}\,(\hat{X}_n\odot_i \hat{X}_m).
\end{equation}
Next, every pair of tree isomorphisms $\T' \cong \T''$, $\sS' \cong \sS''$, leads to an isomorphism $$\T'\odot_i \sS'\cong \T''\odot_i \sS'',$$
i.e.~the image 
\begin{equation*}
\begin{tikzcd}
\it{im}\big(\lT(n) \times \lT(m)  & \lT(m+n-1)\big)
	\arrow["\odot_i",from=1-1, to=1-2]
\end{tikzcd}\end{equation*}
is a full subgroupoid of $\lT(m+n-1)$, and hence there is a map
\begin{equation}\label{eq:composition_in_F_2}
\begin{tikzcd}
\underset{\lT(n)\times\lT(m)}{\it{colim}} (\hat{X}_n \odot_i \hat{X}_m) & \underset{\lT(m+n-1)}{\it{colim}}\hspace{5pt}\hat{X}_{m+n-1}=F_1(X)_{m+n-1}.
	\arrow["\odot_i",from=1-1, to=1-2]
\end{tikzcd}\end{equation}
The composite of (\ref{eq:composition_in_F_1}) and (\ref{eq:composition_in_F_2}) provides the
composition maps $\circ_i$ for the operad $F_1(X)$. It is straightforward to check that the associativity axioms hold and we have:
\begin{lemma}\label{lemma:F1}
The operad $F_1(X)$ with $\circ_i$ defined above is the free non-symmetric non-unital internal operad generated by $X$ in $\OO$, i.e.~$F_1$ is the left adjoint of $U_1$.
\end{lemma}
\pf
The proof is standard. 
Let $A$ be a non-unital non-symmetric internal operad in $\OO$ and $y\colon X \rightarrow A$ a map of internal collections in $\OO$. We have to construct a map of operads $Y\colon F_1(X)\rightarrow A$, such that 
\begin{equation}\label{diagram:free-operad}
\begin{tikzcd}
	& F_1(X) \\
	X & A
	\arrow["i", from=2-1, to=1-2]
	\arrow["y"', from=2-1, to=2-2]
	\arrow["Y", dashed, from=1-2, to=2-2]
\end{tikzcd}
\end{equation}
commutes in the category of internal collections.
We describe the map $Y$ by constructing a map of collections $\hat{X}(\T) \xrightarrow{Y(\T)} A,$ for each tree~$\T$, which will respect the level-interchange isomorphism, hence inducing a map from the colimit $F_1(X)$.

Let $\T \in \lT_n$, $\T =\underline{n}\xrightarrow{(\tau_1,i_1)} \cdots \xrightarrow{(\tau_k,i_k)} \underline{1}$ with fiber sequence
$t_1,\ldots, t_k$, and recall that
$$\hat{X}(\T)= (\cdots\, (X_{t_k}\odot_{i_{k-1}} X_{t_{k-1}})\odot_{i_{k-2}}\cdots )\odot_{i_{1}}X_{t_1}.$$
Let us similarly denote 
$$\hat{A}(\T)= (\cdots\, (A_{t_k}\odot_{i_{k-1}} A_{t_{k-1}})\odot_{i_{k-2}}\cdots )\odot_{i_{1}}A_{t_1}$$
and 
$$\hat{y}(\T)\colon \hat{X}(\T) \rightarrow \hat{A}(\T),$$
$$\hat{y}(\T)= (\cdots\, (y_{t_k}\odot_{i_{k-1}} y_{t_{k-1}})\odot_{i_{k-2}}\cdots )\odot_{i_{1}}y_{t_1}.$$
We can compose the sequence $\hat{A}(\T)$ in the operad $A$ by
$$\hat{A}(\T)\xrightarrow{\circ_{i_{k-1}}\,\cdots\, \circ_{i_1}} A_{n}.$$
Thus, we define $Y(\T)$ as the composite
$$\hat{X}(\T) \xrightarrow{\hat{y}(\T)} \hat{A}(\T)\xrightarrow{\circ_{i_{k-1}}\cdots \circ_{i_1}} A_{n}.
$$
The parallel associativity of the operad $A$ ensures that the diagram
\[\begin{tikzcd}
	{\hat{X}(\T')} & {A_n} \\
	{\hat{X}(\T'')}
	\arrow["{Y(\T')}", from=1-1, to=1-2]
	\arrow["{Y(\T'')}"', from=2-1, to=1-2]
	\arrow["\cong"', from=1-1, to=2-1]
\end{tikzcd}\]
 commutes for any isomorphism $\T'\cong {\T''}$ of $\lT_n$, and hence the maps $Y(\T)$ induce a unique map $Y\colon  F_1(X)\rightarrow A$. By construction, the diagram (\ref{diagram:free-operad}) commutes.
It remains to show that~$Y$ preserves composition. This is ensured by the sequential associativity of the operad $A$.
\epf
The symmetrisation and adding the unit to a non-unital non-symmetric internal operad is analogous to the classical case.
Let $A$ be an internal non-symmetric operad in $\OO$ and define
$$F_2(A)_n = \displaystyle\bigoplus_{\pi\in \Sigma_n}\pi (A_n).$$
The internal actions
\[\begin{tikzcd}
	\sigma (F_2(A)_n)  & F_2(A)_n
	\arrow["\overline{\sigma}", from=1-1, to=1-2]
\end{tikzcd}\]
    are induced by the component maps 
\[\begin{tikzcd}
	\sigma(\pi (A_n)) & \sigma\pi (A_n)  & F_2(A)_n,
	\arrow["\cong", from=1-1, to=1-2]
 \arrow[hook, from=1-2, to=1-3]
\end{tikzcd}\] 
and the composition
\[\begin{tikzcd}
	F_2(A)_n \odot_iF_2(A)_m  & F_2(A)_{m+n-1}
	\arrow["\circ_i", from=1-1, to=1-2]
\end{tikzcd}\]
   is induced by the maps
   \[\begin{tikzcd}
	\sigma(A_n) \odot_i \tau (A_m) & (\sigma\circ_i\tau) (A_n\odot_{\sigma(i)}A_m )  && (\sigma\circ_i\tau)(A_{m+n-1})\\
 &&& F_2(A)_{m+n-1}.
	\arrow["\cong", from=1-1, to=1-2]
 \arrow["(\sigma\circ_i\tau)(\circ_{\sigma(i)})", from=1-2, to=1-4]
  \arrow[hook, from=1-4, to=2-4]
  \arrow[dashed,curve={height=25pt}, from=1-1, to=2-4]
\end{tikzcd}\] 
\begin{lemma}\label{lemma:F2}
    The functor $F_2$ is left adjoint to the forgetful functor $U_2$.
\end{lemma}
\pf
The verification is straightforward.
\epf
For the categorical operad $\SSS$, Lemma~\ref{lemma:F2} recovers the symmetrisation of~\cite{Day:substitution}.
Let~$U$ denote the unit of a unital categorical operad $\OO$. We define an internal unital operad~$I$ in~$\OO$ by $$I_1:=U,$$ and for $n\neq1$, $$I_n:= \emptyset,$$ the initial object of $\OO(n)$.
\begin{lemma}\label{lemma:F3}
The coproduct $F_3(A) = A\oplus I$ defines left adjoint to the forgetful functor~$U_3$.
\end{lemma}
\pf
    The verification is straightforward.
\epf

Putting together Lemmas \ref{lemma:F1}, \ref{lemma:F2}, and \ref{lemma:F3} we obtain:
\begin{proposition}
The composite $F= F_3\circ F_2\circ F_1$ is left adjoint to the functor $$U\colon \it{IOp}(\OO)\rightarrow \it{IColl}(\OO).$$
\end{proposition}
\begin{example}\label{example:free_operad}
We analyze the free construction for $\OO=\SSS$, the operad defined in Section~\ref{section:SSS}, and $\V=\Vect$, the category of vector spaces over a field $\mathbb{k}$. This gives a description of the free $C$-operad over a non-symmetric $C$-collection $X$. Fix a leveled tree $$\T =\underline{n}\xrightarrow{(\tau_1,i_1)} \cdots \xrightarrow{(\tau_{k},i_{k})} \underline{1}$$ with fiber sequence
$t_1,\ldots, t_k.$
The tree $\T$ can be regarded as a planar graph with an order of vertices given by the levels of $\T$. We label each edge of $\T$ by an object of $C$, i.e.~we fix a labeling function 
$$\lambda\colon \it{Edges}(\T) \rightarrow \it{ob}(C).$$ 
Then each vertex $v_j$ has a type $$\overline{v}_j=\colorop (c_1\cdots c_{t_j};c)$$ given by the labels of adjacent edges.
Let $X$ be a non-symmetric $C$-collection of functors 
\[\begin{tikzcd}
\underbrace{C^{\it{op}} \times \cdots \times C^{\it{op}}}_{n-times} \times C & \it{Vect}_{\mathbb{k}}.
\arrow["X_n", from=1-1, to=1-2]
\end{tikzcd}\]
We define
$$\hat{X}(\T,\lambda) :=X_{t_k}(\overline{v}_k)\otimes \cdots \otimes X_{t_1}(\overline{v}_1).$$
The value $F_1(X)\colorop (l_1\,\cdots \,l_n;r)$ is the coproduct $$F_1(X)\colorop (l_1\,\cdots \,l_n;r)=\displaystyle\bigoplus_{(\T,\lambda)}\hat{X}(\T,\lambda), $$ taken over all leveled planar trees with $n$ leaves, such that $\lambda$ labels the leaves of $\T$ by the objects $l_1,\ldots, l_n\in C$ and the root edge by $r\in C$, with the following identifications:
\begin{itemize}
    \item Let $\T'\cong \T''$, such that the two trees differ only in the $i$-th position by a square of elementary morphisms with disjoint fibers as in (\ref{isomorphic-trees}). Then 
    $$ a_k\otimes \cdots \otimes a_i \otimes a_{i-1} \otimes \cdots \otimes a_1 \in \hat{X}(\T',\lambda) $$
    is identified with
    $$ a_k\otimes \cdots \otimes a_{i-1} \otimes a_{i} \otimes \cdots \otimes a_1 \in \hat{X}(\T'',\lambda). $$
    
    \item Let an edge $e$ of a labeled tree $(\T,\lambda)$ connect vertices $v_i$ and $v_j$ (with $v_j$ closer to the root) of types $$\overline{v}_j=\colorop (y_1\,\cdots \, d \,\cdots\,y_{t_j};y),\overline{v}_i=\colorop (x_1\,\cdots \,x_{t_i};c) $$
    and let $f\colon c \rightarrow d$ in $C$. Then 
    $$ a_k\otimes \cdots \otimes a_j\cdot f \otimes\cdots \otimes  a_i \otimes \cdots \otimes a_1 $$
    is identified with
    $$ a_k\otimes \cdots \otimes a_{j} \otimes \cdots\otimes f \cdot a_{i} \otimes \cdots \otimes a_1. $$
\end{itemize}
The notation $a_j\cdot f$ and $ f\cdot a_{i}$ was introduced in Example \ref{example:explicit_odot}. Note, that in the graded context signs may be involved.
The free symmetric $C$-operad $F(X)$ has $$F(X)\colorop (l_1\,\cdots \,l_n;r)=\displaystyle\bigoplus_{\sigma \in \Sigma_n}F_1(X)\colorop (l_{\sigma(1)}\,\cdots \,l_{\sigma(n)};r).$$
Elements of $F(X)$ are thus sums of trees as in Figure~\ref{fig:labeled_tree}, with specified orders on their leaves.
\end{example}

\section{Algebras}\label{section:algebras}
The first part of this section modifies some standard notions to the category-colored setting. In the second part we express algebras of classical operads with values in a functor category~$\V^C$ as algebras of $C$-operads in $\V$. For this we need $\V$ to be a cocomplete symmetric closed monoidal category  and $C$ a small symmetric monoidal $\V$-category, so that the category $\V^C$ is symmetric monoidal closed via the Day convolution \cite{Day:closed}.

From this point on we work only with the monoidal category $\V=\it{Vect}_{\mathbb{k}}$, the category of vector spaces over a field ${\mathbb{k}}$. This allows us to work with elements, which simplifies the presentation. We believe that the reader can easily translate the definitions and constructions to other relevant monoidal categories.

The category of algebras of a $C$-operad $P$ lives over the functor category $\it{Vect}_{\mathbb{k}}^C.$ For a functor $A\colon C \rightarrow \Vect$, we define its \textit{endomorphism $C$-operad} $\it{End}_A$. Let $c_1, \ldots, c_n, c$ be objects of $C$ and let $\it{Hom}(U,V)$ denote the space of linear maps from $U$ to $V$. We define 
	$$\colorop \it{End}_A(c_1\, \cdots\, c_n ;c) := \it{Hom}(A(c_1)\otimes\cdot\cdot\cdot\otimes A(c_n), A(c)).$$
For a map 
	\[\begin{tikzcd}
\colorop (\sigma\, f_1 \,\cdots \,f_n;f)\colon \colorop (c_1  \,\cdots \,c_n;c)
& \colorop (d_1 \,\cdots\, d_n;d)
	\arrow[from=1-1, to=1-2]
\end{tikzcd}\]
in $(\bS C)^{\it{op}}\times C$ and $\varphi \in \colorop \it{End}_A(c_1\, \cdots\, c_n ;c)$, we define $ \colorop \it{End}_A(\sigma \, f_1 \cdots f_n;f)(\varphi) $ to be the composite
\[\begin{tikzcd}
	{A(c_1)\otimes\cdots\otimes A(c_n)} &&& {A(c)} \\
 {A(d_{\sigma(1)})\otimes\cdots\otimes A(d_{\sigma(n)})} &&& \\
	{A(d_1)\otimes\cdots\otimes A(d_n)} &&& {A(d)} .
	\arrow[" \colorop \it{End(A)}(\sigma\,f_1 \cdots f_n;f)(\varphi)"', dashed, from=3-1, to=3-4]
	\arrow["\varphi", from=1-1, to=1-4]
	\arrow["{A(f)}",from=1-4, to=3-4]
	\arrow["(-)\sigma",from=3-1, to=2-1]
	\arrow["{A(f_1)\otimes\cdots\otimes A(f_n)}",from=2-1, to=1-1]
\end{tikzcd}\]
	For $n=0$, let $ \colorop \it{End}_A(\emptyset;c):= A(c)$.
Let $\varphi \in \colorop \it{End}_A(c_1\cdots c_n;c)$ and $\psi \in \colorop
\it{End}_A(d_1\cdots d_m;c_i)$. The operad composition is defined by
$$\varphi\circ_i \psi  := \varphi\circ(\uu_{A(c_1)}\otimes \cdots \otimes \psi \otimes
\cdots \otimes \uu_{A(c_n)}).$$
The operad $\it{End}_A$ is unital with units $A(f) \in \colorop \it{End}_A(c;d)$, for $f\colon c\to d$ in $C$.
\begin{definition}
	Let $P$ be a $C$-operad.
	A $P$-algebra is a functor $A\colon C \rightarrow \Vect$ together with a morphism of $C$-operads 
 \[\begin{tikzcd}
P  & \it{End}_A.	\arrow["\alpha",from=1-1, to=1-2]
\end{tikzcd}\]
\end{definition}
A $P$-algebra is thus a functor $A\colon C \rightarrow \Vect$, together with maps 
$$
\displaystyle\int^{c_1\cdots c_n \in \bS C^{\it{op}}}\colorop P(c_1 \, \cdots \, c_n;c)\otimes A(c_1)\otimes \cdots \otimes A(c_n) \xrightarrow{{\alpha}}A(c),
$$
for any $n\geq0$ and object $c \in C$, which satisfy some obvious axioms.
Due to our Proposition~\ref{prop:C-operads_and_Substitudes}, Proposition 5.1.6 of \cite{Batanin_White:substitudes} offers other characterisations of algebras of unital $C$-operads.



\begin{definition}
Let $F\colon D \rightarrow \Vect$ be a functor.
An equivalence relation $E$ on the functor $F$ is a subfunctor $$E(d)\subseteq F(d)\times F(d),$$
such that, for every $d\in D$, $E(d)$ is an equivalence relation.
The quotient $F/E$ is the functor $$F/E(d):=F(d)/E(d).$$
\end{definition}
\begin{definition}
 Let $R$ be a set of pairs 
$$R=\{(a_s,b_s)\,|\,a_s,b_s \in F(c_s)\}_{s\in S},$$
for some set $S$.
An equivalence relation generated by $R$ is the smallest equivalence relation on $F$ which contains $R$. The generating pairs will be denoted $a_s\sim b_s$.
\end{definition}
For a non-symmetric $C$-collection $X$, the notation $\alpha\in X$ will mean $\alpha \in \colorop X_n(c_1\, \cdots\, c_n ;c)$, for some $n\in \NN$ and objects $c_1, \ldots , c_n , c$ of $C$.
\begin{definition}
    An operadic ideal $I$ of a $C$-operad $P$ is a subfunctor $I\subseteq P$, 
    such that every composition
$\gamma\circ_i\beta$ of two elements of $P$, where at least one of the elements is in $I$, is again in $I$.\end{definition}
   
\begin{definition}
An operadic ideal of a $C$-operad $P$ generated by a set
$$
A=\{\alpha_s \in P\}_{s\in S},
$$ for some set $S$,
is the smallest operadic ideal of $P$, which contains $A$. We will denote it by~$(A)$.
\end{definition}
For an operadic ideal $I$ of a $C$-operad $P$, the $C$-collection $P/I$ acquires an obvious \hbox{$C$-operad} structure. For a non-symmetric $C$-collection $X$ in~$\Vect$, every operation \hbox{$\alpha\in F(X)$} is a finite sum $\alpha=\sum_{s\in S}\alpha_s$ with each $\alpha_s\in\hat{X}(\T_s, \lambda_s)$ for some tree $\T_s$ with labeling $\lambda_s$, cf.~Example~\ref{example:free_operad}.
\begin{definition}
We say that $\alpha\in F(X)$ is homogeneous of weight $k$, if every tree $\T_s$ appearing in the sum $\alpha=\sum_{s\in S}\alpha_s$, with $\alpha_s\in\hat{X}(\T_s, \lambda_s),$ has $k$ vertices. We denote by~$F(X)^{(k)}$ the subfunctor of $F(X)$ of operations of weight~$k$.
\end{definition}
\begin{definition}
    A $C$-operad $P$ is quadratic binary if it is of the form $P=F(X)/(A)$, such that $X_n$ is trivial unless $n=2$, and $A\subseteq F(X)^{(2)}$. 
\end{definition}
Let $(C,\oplus,0)$ be a small symmetric monoidal linear category (i.e.~a $\Vect$-enriched category), and $P$ an operad with values in the category $\it{Vect}_{\mathbb{k}}^C$ of linear functors and their linear natural transformations.
With the Day convolution, $\it{Vect}_{\mathbb{k}}^C$ is a cocomplete symmetric monoidal closed category.
Let $P(n,r)$ denote the value of $P(n)$ on an object~$r$ of~$C$.
Our goal is to define a $C$-operad $P^C$ with values in $\Vect$ which has the same algebras as the operad $P$.
The object $\colorop P^C(c_1 \, \cdots \, c_n;c)$ is the linear coend
$$
\displaystyle\int^{r\in C}C(c_1\oplus\cdots \oplus
c_n\oplus r,c)\otimes P(n,r),
$$
and the action of a morphism $\colorop (\sigma \,\alpha_1 \, \cdots \, \alpha_n;\hspace{8pt}\alpha)$ of $\bS C^{\it{op}}\times C$ is induced by precomposition with 
$((\alpha_1 \oplus \cdots \oplus \alpha_n) \circ \overline{\sigma}) \oplus 1_r$ and postcomposition with $\alpha$ on  $C(c_1\oplus\cdots \oplus c_n\oplus r,c)$.
The composition maps $\circ_i$ of $P^C$ are induced by the composition maps $\circ_i^P$ of $P$, the composition $\circ^C$ of $C$, and the symmetry isomorphism $s$ of $\V$:

\[\begin{tikzcd}
	{C(c_1\oplus\cdots \oplus c_n\oplus p,c)\otimes P(n,p) \otimes C(d_1\oplus\cdots \oplus d_m\oplus q,c_i)\otimes P(m,q)} \\
	{C(c_1\oplus\cdots \oplus c_n\oplus p,c)\otimes C(d_1\oplus\cdots \oplus d_m\oplus q,c_i) \otimes P(n,p)\otimes P(m,q)} \\
	{C(c_1\oplus\cdots\oplus d_1\oplus\cdots \oplus d_m\oplus\cdots\oplus c_n\oplus p\oplus q,c)\otimes P(n + m-1,p\oplus q)}.
	\arrow["{\uu\otimes s \otimes \uu}", from=1-1, to=2-1]
	\arrow["{\circ^C\otimes \circ_i^P}", from=2-1, to=3-1]
\end{tikzcd}\]

\begin{proposition}\label{prop:algebras_as_algebras}
    The categories of $P$-algebras in $\it{Vect}_{\mathbb{k}}^C$ and $P^C$-algebras in $\Vect$ are isomorphic.
\end{proposition}
\pf
    
By definition, a $P$-algebra is a linear functor $X$ together with a linear natural transformation 
$$
\displaystyle\int^{r\in C, c_1\cdots c_n \in \bS C^{\it{op}}}C(c_1\oplus\cdots \oplus
c_n\oplus r,c)\otimes P(n,r)\otimes X(c_1)\otimes \cdots \otimes X(c_n) \xrightarrow{\alpha_c}X(c),
$$
natural in $c$, which preserves the operad structure of $P$.
This is the same as
$$
\displaystyle\int^{c_1\cdots c_n \in \bS C^{\it{op}}}\colorop P^C(c_1 \, \cdots \, c_n;c)\otimes X(c_1)\otimes \cdots \otimes X(c_n) \xrightarrow{{\alpha_c}}X(c).
$$
The map ${\alpha_c}$ respects the $C$-operad structure, since the compositions of $P^C$ are induced by the composition of $P$, and hence $(X,\alpha)$ is a $P^C$-algebra. This correspondence induces isomorphism of the categories of algebras.
\epf
\begin{example}\label{example:dgas}
As an instance of Proposition~\ref{prop:algebras_as_algebras}, we show that differential graded associative algebras are algebras of a non-unital $\Vect$-enriched category-colored operad. The coloring category in this case is not a groupoid, which is related to the question of Remark~3.12 of~\cite{petersen2013operad}.
Let us define the coloring linear category $D$. Objects of $D$ are integers, the spaces $D(n,n-1)$, for any integer $n$, are one dimensional with generator~$\partial_n$, and the other Hom-spaces of $D$ are trivial. This forces $\partial\partial=0$, hence linear functors $D\rightarrow\Vect$ are exactly chain complexes of vector spaces.
We define a linear functor 
$$\displaystyle\bigoplus_{n\in \NN}(D^{\it{op}})^n\times D \xrightarrow{X}\Vect,$$
\begin{center}
\begin{tabular}{ c c }
 $\colorop X(m\hspace{8pt}n;m+n):= \mathbb{k}[\mu_{m,n}]$,&
$\colorop X(m\hspace{8pt}n;m+n-1):= \mathbb{k}[\mu_0,\mu_1,\mu_2]/\langle\mu_1+(-1)^{n}\mu_2-\mu_0\rangle,$
\end{tabular}
\end{center}
for any $m,n\in D$,
and $X$ is trivial otherwise.
The only non-trivial $D$-actions are
\begin{center}
\begin{tabular}{ c c c }
$\colorop X(\partial\,\uu;\uu)(\mu_{m-1,n}):=\mu_1,$&
$\colorop X(\uu\,\partial;\uu)(\mu_{m,n-1}):=\mu_2,$&
$\colorop X(\uu\,\uu;\partial)(\mu_{m,n}):=\mu_0.$
\end{tabular}
\end{center}
The functor $X$ generates the free non-unital $D$-operad $F(X)$ which we further divide by an operadic ideal $I$ generated by the elements
$$
\mu_{m+n,k}\circ_1\mu_{m,n} - \mu_{m,n+k}\circ_2 \mu_{n,k},
$$
for all $m,n,k \in D$.
One can check, that algebras of the $D$-operad $F(X)/I$ are precisely differential graded associative algebras.
\end{example}

\section{Hyperoperad for operads}\label{section:operad_for_operads}
We construct a quadratic operad $\bH$, whose algebras are classical symmetric non-unital Markl operads. This serves as a toy example. By the same procedure, we describe a quadratic operad $\bH_{\bO}$, whose algebras are non-unital Markl operads in the operadic category $\bO$, and further a quadratic operad $\bH_C$, for a fixed small category $C$, whose algebras are (non-unital) $C$-operads. We will work with the category \hbox{$\V=\Vect$}, believing that the reader can easily reformulate the construction into other monoidal categories.  
\subsection*{Classical non-unital Markl operads}
\begin{figure}
    \centering
    \includegraphics{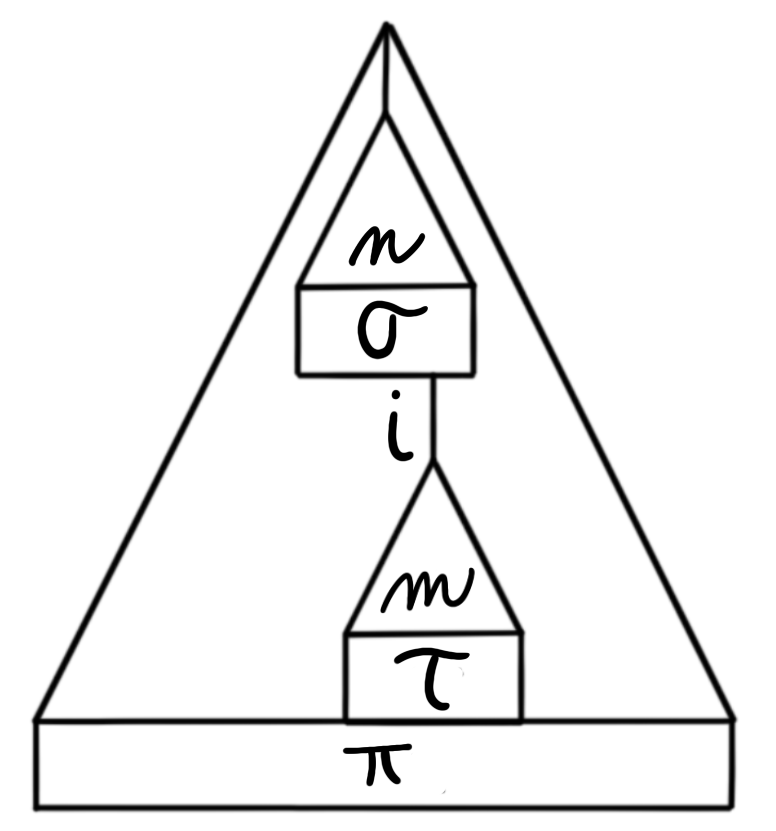}
    \caption{Generators of the operad $\bH$.}
    \label{fig:generators}
\end{figure}

\begin{figure}
    \centering
    \includegraphics{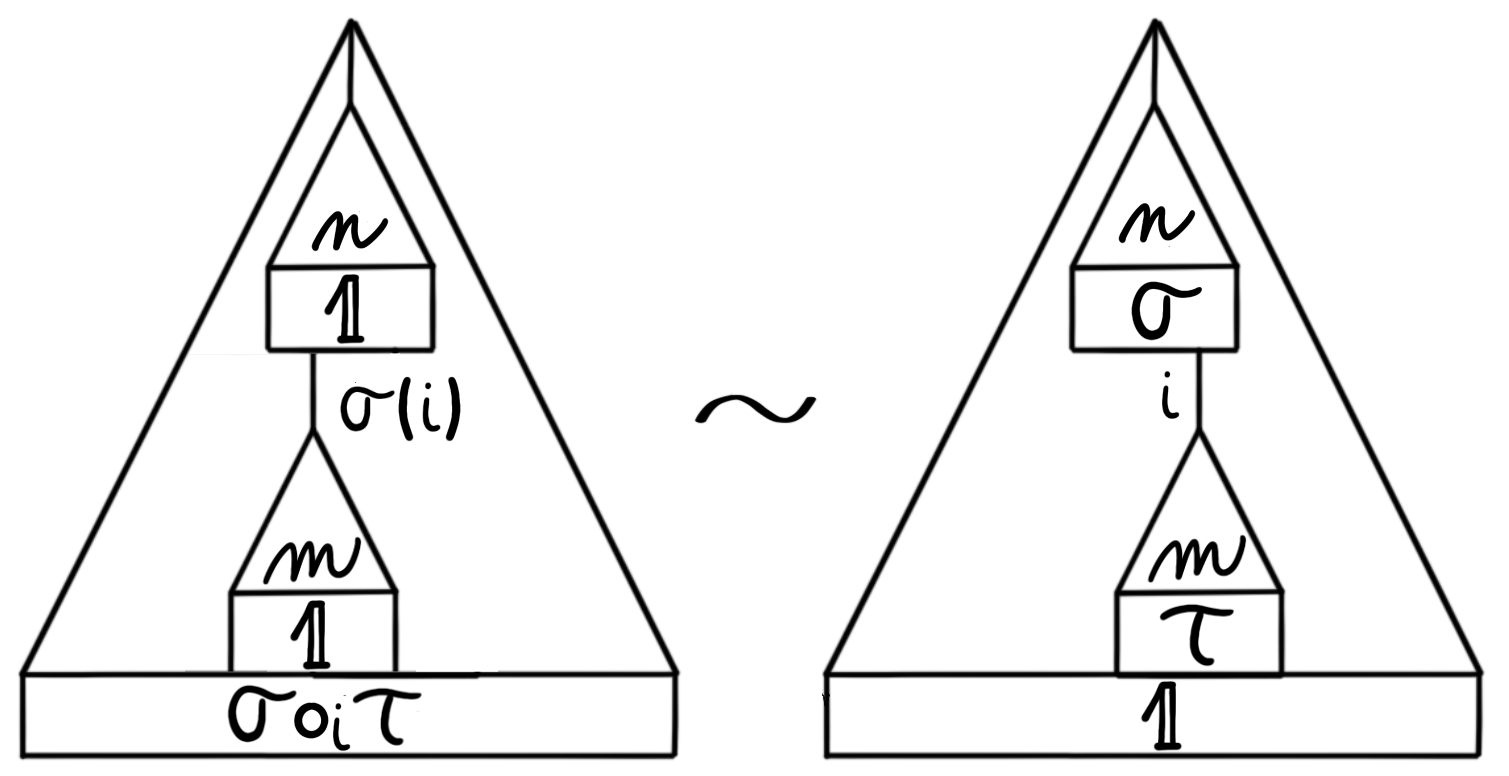}
 \caption{Equivariance of the symbols $*_i$.}
    \label{fig:equiavariance_of_*_i}
\end{figure}

\begin{figure}
    \centering
    \includegraphics{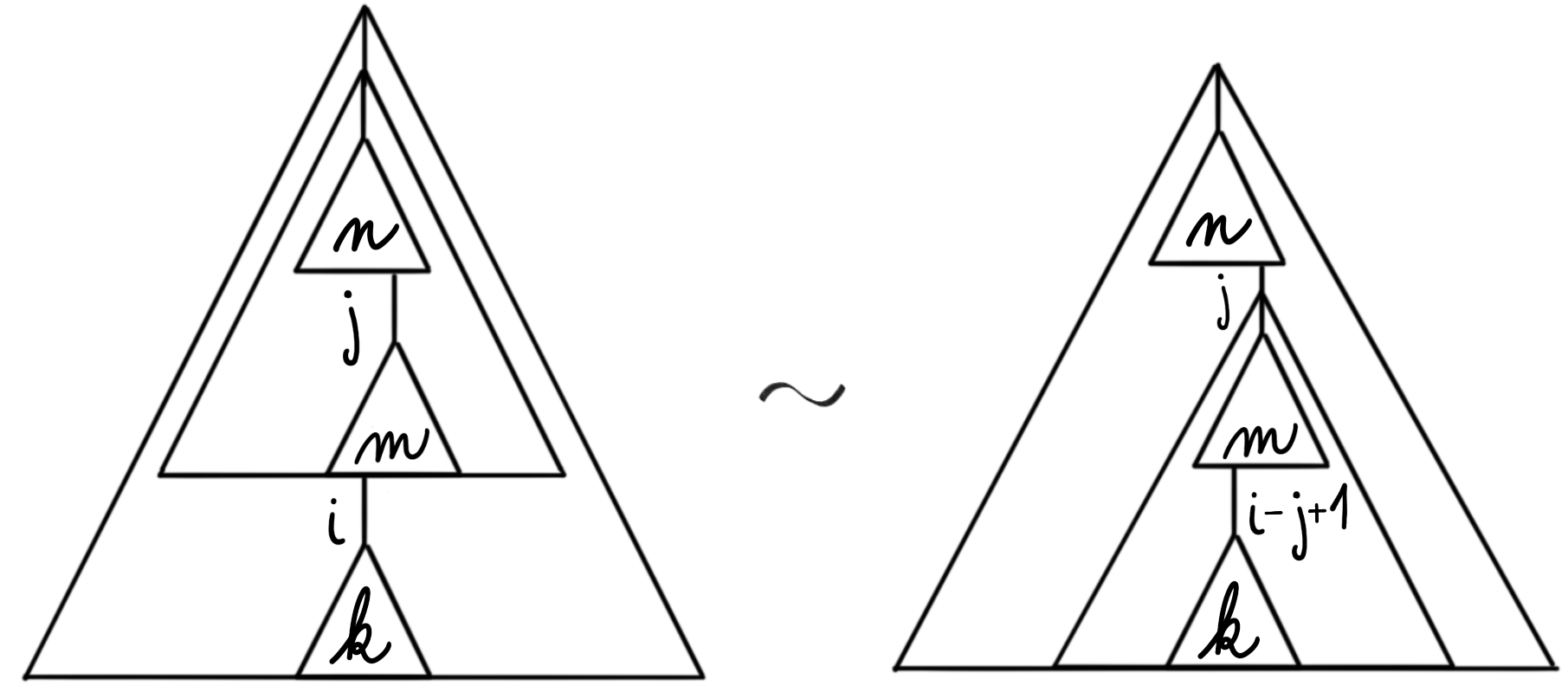}
    \caption{Parallel associativity identification.}
    \label{fig:identification_seq}
\end{figure}

We denote by~$\Sigma$ the free symmetric monoidal category on one generator. Its objects are natural numbers and morphisms are permutations. This will be the category of colors for the operad $\bH$.
Let us begin with a generating non-symmetric $\Sigma$-collection 
\[
\begin{tikzcd}
	(\Sigma^{\it{op}})^{k} \times \Sigma & \Vect,
	\arrow["X_k", from=1-1, to=1-2]
\end{tikzcd}
\]
defined by:
\begin{itemize}
    \item For any $m\geq 0, n\geq 1$, the space 
$\colorop X_2(n\hspace{.5cm}m;n+m-1)$ is freely generated by the set
$$\{*_i\, |\, 1\leq
i\leq n\} \times \Sigma_n\times \Sigma_m\times\Sigma_{n+m-1}.$$
The generators will be written as schemes 
$\colorop (*_i,\sigma\,\tau;\hspace{15pt}\pi)$ and the actions of $\Sigma$ on individual inputs and output are given by
$$\colorop X_2(\sigma'\tau';\pi') \colorop (*_i,\sigma\,\tau;\hspace{15pt}\pi) := \colorop (*_i,\sigma\sigma'\,\tau\tau';\hspace{15pt}\pi'\pi).$$
 \item $X_k$ is trivial for $k\neq 2$.
\end{itemize}
Each generator $\colorop (*_i,\sigma\,\tau;\hspace{15pt}\pi)$ represents the composite 
\[
\begin{tikzcd}
	 P(n)\otimes P(m) & P(n+m-1)\\ P(n)\otimes P(m)  & P(n+m-1)
	\arrow["(-)\sigma \otimes (-)\tau", from=2-1, to=1-1]
 \arrow["\circ_i", from=1-1, to=1-2]
 \arrow["(-)\pi", from=1-2, to=2-2]
 \arrow[dashed, from=2-1, to=2-2]
\end{tikzcd}
\]
of symmetric actions and the operation $\circ_i$ of a classical Markl operad $P$. The generators are depicted in Figure~\ref{fig:generators}. 
Next, we form the quotient $Q=X/\textit{Eq}$, where
$\textit{Eq}$ is the equivalence relation, generated by pairs
$$\colorop (*_i,\sigma\,\tau;\hspace{15pt}\uu)\sim \colorop(*_{\sigma(i)},\uu\hspace{10pt}\uu;\hspace{28pt}\sigma\circ_i\tau),$$
for each $m\geq0, 1\geq i\geq n, \sigma \in \Sigma_n$ and $\tau \in \Sigma_m$ (cf.~Figure~\ref{fig:equiavariance_of_*_i}). The permutation $\sigma\circ_i\tau$ is given by inserting the permutation $\tau$ at the $i$-th place in $\sigma$.
Each space $\colorop Q_2(n\hspace{.5cm}m;n+m-1)$ is thus generated by equivalence
classes $\coloropsq [*_i, \sigma\,\tau;\hspace{15pt}\pi]$. Let us denote $[*_i]$~$:=$~$\coloropsq [*_i, \uu\,\uu;\hspace{15pt}\uu].$
Consider the free symmetric $\Sigma$-operad $F(Q)$ (cf.~Example~\ref{example:free_operad}) and form an operadic ideal $\textit{As}$ generated by the elements
\begin{equation}\label{eq:associativity_of_*}
\begin{aligned}
  [*_i]\circ_1[*_j]&-\tau([*_{j+k-1}]\circ_1[*_i]), 
    \\
    [*_i] \circ_1 [*_j]&-[*_j] \circ_2 [*_{i-j+1}],
    \\
    [*_i]\circ_1[*_j]&-\tau([*_j]\circ_1[*_{i-m+1}]),
\end{aligned}
\end{equation}
for all $i,j,k,m,n$ as in (\ref{eq:associativity}).
The symbol $\tau$ stands for the map 
\[
\begin{tikzcd}
	\colorop{F(Q)}(n\hspace{.5cm}m\hspace{.5cm}k;n+m+k-2) &&&& \colorop {F(Q)}(n\hspace{.5cm}k\hspace{.5cm}m;n+m+k-2),
	\arrow["{\colorop {F(Q)}(\tau\,\uu\,\uu\,\uu;\hspace{8pt}\uu)}", from=1-1, to=1-5]
\end{tikzcd}
\]
which swaps the second and third input. An element of the second type is depicted in~Figure~\ref{fig:identification_seq}.

\begin{theoremB1}\label{theorem:operads_as_algebras}
 The category of algebras of the binary quadratic $\Sigma$-operad $$\mathbb{H}:= F(Q)/\textit{As}$$ is isomorphic to the category of classical symmetric non-unital Markl operads.
\end{theoremB1}	
 \pf
The theorem follows directly from definitions and the construction of $\mathbb{H}$.
The~data of an $\mathbb{H}$-algebra $(M,\alpha)$ consists of a symmetric collection of vector spaces
\[
\begin{tikzcd}
	\Sigma & \Vect
	\arrow["M", from=1-1, to=1-2]
\end{tikzcd}
\]
 and of an operad map 
 \[
\begin{tikzcd}
	  \bH & \Endm.
	\arrow["\alpha", from=1-1, to=1-2]
\end{tikzcd}
\]
The map $\alpha$ interprets each class $[*_i] \in \colorop \bH(n\hspace{.5cm}m;n+m-1)$ as a map $\circ_i:=\alpha([*_i])$,
\[
\begin{tikzcd}
	 M(n)\otimes M(m) &  M(n+m-1).
	\arrow["\circ_i", from=1-1, to=1-2]
\end{tikzcd}
\]
Since $\alpha$ is a map of $\Sigma$-collections, we have 
 \begin{align*}
     \circ_i \cdot \big((-)\sigma\otimes (-)\tau\big)
     &=
\Endm\colorop(\sigma\,\tau;\uu)(\circ_i)
=\big(\Endm\colorop(\sigma\,\tau;\uu)\cdot \alpha\big)([*_i])\\
 &=
 \big(\alpha\cdot \colorop \bH(\sigma\,\tau;\uu)\big)\big(\coloropsq [*_i,\uu\,\uu;\hspace{15pt}\uu]\big)
 =
  \alpha\big(\coloropsq [*_i,\sigma\,\tau;\hspace{15pt}\uu]\big)
 =
  \alpha\big( \coloropsq [*_{\sigma(i)},\uu\hspace{10pt}\uu;\hspace{28pt}\sigma\circ_i\tau]\big)
 \\
 &=
 \big(\alpha\cdot \colorop \bH(\uu\hspace{10pt}\uu;\sigma\circ_i\tau)\big)\big(\coloropsq [*_{\sigma(i)},\uu\,\uu;\hspace{28pt}\uu]\big)
 =
 \Endm\colorop(\uu\hspace{10pt}\uu;\sigma\circ_i\tau)(\circ_{\sigma(i)})\\
 &=
\big((-)(\sigma\circ_i\tau)\big) \cdot \circ_{\sigma(i)}.
 \end{align*}
 
This shows that the operations $\circ_i$ are equivariant.
The associativity of the operations $\circ_i$ is ensured by the identifications (\ref{eq:associativity_of_*}) and the fact, that $\alpha$ is an operad morphism. On~the other hand, one can reconstruct the map $\alpha$ from a non-unital operad structure on~$M$. Further, the morphisms of non-unital operads and of $\mathbb{H}$-algebras agree.
 \epf

\subsection*{Non-unital Markl $\bO$-operads}
The construction in the previous section generalises to the context of non-unital Markl operads in operadic categories. We use freely the terminology of \cite{Batanin_Markl:kodu2022,Batanin_Markl:kodu2021}.

Let $\bO$ denote an operadic category, which satisfies the axioms required for the definition of Markl $\bO$-operads \cite[Section~6]{Batanin_Markl:kodu2022}. That is, we assume that $\bO$ is a graded, factorizable operadic category in which all quasibijections are invertible, the strong blow-up and
unique ﬁber axioms are fulfilled, and a morphism $f$ is an isomorphism if and only it is of grade~0. We also assume that $\bO$ is small, but this is the case for the examples of \cite{Batanin_Markl:kodu2022,Batanin_Markl:kodu2021}. Let~$\bO_{\it{iso}}$ denote the full subcategory of isomorphisms of $\bO$. An elementary morphism $\phi\colon T \rightarrow S$ with the non-trivial fiber $F$ over $i\in |S|$ will be written as 
$
\phi\colon F \triangleright_i T \rightarrow S$ or simply $\phi\colon F \triangleright T \rightarrow S
$, when the concrete $i$ is not important.
We construct an $\bO_{\it{iso}}^{\it{op}}$-operad $\bH_\bO$ whose algebras are (non-unital) Markl $\bO$-operads.

Let $Bq(\bO_{\it{iso}}^{\it{op}})$ denote the category $\bigoplus_{n\geq0} (\bO_{\it{iso}})^{\times n} \times \bO_{\it{iso}}^{\it{op}}$ and $Bq^{\textit{Disc}}(\bO_{\it{iso}}^{\it{op}})$ the discrete category of objects of $Bq(\bO_{\it{iso}}^{\it{op}})$. We define a collection 
$X^{\textit{Disc}}\colon Bq^{\textit{Disc}}(\bO_{\it{iso}}^{\it{op}}) \rightarrow \Vect$
by
\begin{itemize}
    \item 
 $\colorop X^{\textit{Disc}}(S\,F;T)$ is generated by the set $\{*_{\phi} \,|\, \phi\colon F \triangleright T \rightarrow S \text{ is
elementary in $\bO$}\}$,
\item $X^{\textit{Disc}}$ is trivial otherwise.
\end{itemize}
We equip the generating collection with free $\bO_{\it{iso}}^{\it{op}}$-actions on inputs and output, which is done using Kan extensions. The left Kan extension $X=\textit{Lan}(X^{\textit{Disc}})$ of the functor $X^{\textit{Disc}}$ along the inclusion functor $$Bq^{\textit{Disc}}(\bO_{\it{iso}}^{\it{op}})
\subseteq Bq(\bO_{\it{iso}}^{\it{op}})$$ is then a non-symmetric $\bO_{\it{iso}}^{\it{op}}$-collection. Explicitly, $X$ is given by
\begin{align}\label{equation:Kan_Extension}
     \colorop X(S' \,F';T')= \big\{\colorop (*_{\phi},\sigma\,\tau;\hspace{15pt}\pi)\,|\,& \tau \colon F' \rightarrow F'', \sigma\colon S' \rightarrow S'', \pi\colon T''\rightarrow T'\text{ are isomorphisms of $\bO$,}\\
& \phi\colon F'' \triangleright T'' \rightarrow S'' \text{ is
elementary in $\bO$}\big\}, \nonumber
\end{align}
and the $\bO_{\it{iso}}$-actions are 
$$\colorop X(\sigma'\,\tau';\pi')\colorop (*_{\phi},\sigma\,\tau;\hspace{15pt}\pi)=\colorop (*_{\phi},\sigma \sigma'\,\tau \tau';\hspace{15pt}\pi' \pi ).$$

Recall from \cite[Definition 6.1]{Batanin_Markl:kodu2022} that for any square

\[
\begin{tikzcd}
	T' & {T''} \\
	S' & {S''}
	\arrow["\phi'", from=1-1, to=2-1]
	\arrow["{\phi''}", from=1-2, to=2-2]
	\arrow["\cong"', from=1-1, to=1-2]
	\arrow["\omega", from=1-1, to=1-2]
	\arrow["\sigma", from=2-1, to=2-2]
	\arrow["\sim"', from=2-1, to=2-2]
\end{tikzcd}
\]
in $\bO$ with $\omega$ an isomorphism, $\sigma$ a quasibijection, and $\phi',\phi''$ elementary, there is an induced isomorphism $\overline{\omega}$
between the unique non-trivial fibers $F'$ and $F''$ of $\phi'$~and~$\phi''$. 
We form the quotient $Q=X/\textit{Eq}$, where $\textit{Eq}$ is an equivalence relation generated by \begin{equation}\label{eq:equivariance_of_*_operadic}
\colorop (*_{\phi'},\uu\,\uu;\hspace{18pt}\omega)
\sim
\colorop (*_{\phi''},\sigma\,\,\overline{\omega};\hspace{15pt}\uu),
\end{equation} 
for any $\sigma,\omega,\phi',\phi''$ as above.
Denote the class of $\colorop (*_{\phi},\uu\,\uu;\hspace{15pt}\uu)$ by $[*_{\phi}]$.
We consider the free symmetric $\bO_{\it{iso}}$-operad $F(Q)$, which we divide by the operadic ideal generated by two types of associators:

\textit{Sequential associators.}
Let $\phi\colon F\triangleright_j T
\rightarrow S$ and $\psi\colon G\triangleright_i S \rightarrow P$ be two elementary maps in~$\bO$, such that  $|\psi|(j)=i$. Their composite is an elementary map
$$\xi=\psi\phi\colon H\triangleright_i T \rightarrow P$$ and there is an induced elementary map $\phi_i\colon F\triangleright H \rightarrow G$: 
\[\begin{tikzcd}
	F & H & G \\
	F & T & S \\
	&& P.
	\arrow["\phi", from=2-2, to=2-3]
	\arrow["\psi", from=2-3, to=3-3]
	\arrow["\psi\phi"', curve={height=18pt}, from=2-2, to=3-3]
	\arrow["\triangleright"{marking}, draw=none, from=1-2, to=2-2]
	\arrow["\triangleright"{marking}, draw=none, from=1-3, to=2-3]
	\arrow["\triangleright"{marking}, draw=none, from=2-1, to=2-2]
	\arrow["{\phi_i}", from=1-2, to=1-3]
	\arrow[no head, from=1-1, to=2-1]
	\arrow[shift left=1, no head, from=1-1, to=2-1]
	\arrow["\triangleright"{marking}, draw=none, from=1-1, to=1-2]
\end{tikzcd}\]
We have
\begin{center}
\begin{tabular}{ c c c c}
 $[*_{\phi}] \in \colorop Q(S\,F;T),$&
$[*_{\psi}] \in \colorop Q(P\,G;S),$ &
$[*_{\phi_i}] \in \colorop Q(G\,F;H),$&
$[*_{\xi}] \in \colorop Q(P\,H;T),$
\end{tabular}
\end{center}
and the composites $[*_{\phi}] \circ_2 [*_{\psi}]$ and $[*_{\xi}] \circ_1
[*_{\phi_i}]$ are both elements of ${F(Q)}\colorop (P\,G\,F;T).$
The sequential associators are
\begin{equation}\label{equation:ideal1}
[*_{\phi}] \circ_2 [*_{\psi}]-[*_{\xi}] \circ_1 [*_{\phi_i}],\end{equation}
for any $\phi$ and $\psi$ as above.

\textit{Parallel associators.} Let $\psi',\psi'',\phi',$ and $\phi''$ form a square of elementary morphisms with disjoint fibers:
\[\begin{tikzcd}
	& {F''} & {G'} \\
	{F'} & T & {P'} \\
	{G''} & {P''} & S.
	\arrow["{\phi'}", from=2-2, to=2-3]
	\arrow["{\psi'}", from=2-3, to=3-3]
	\arrow["{\psi''}"', from=3-2, to=3-3]
	\arrow["{\phi''}"', from=2-2, to=3-2]
	\arrow["\triangleright"{marking}, draw=none, from=1-3, to=2-3]
	\arrow["\triangleright"{marking}, draw=none, from=2-1, to=2-2]
	\arrow["\triangleright"{marking}, draw=none, from=1-2, to=2-2]
	\arrow["\triangleright"{marking}, draw=none, from=3-1, to=3-2]
	\arrow[shift right=1, no head, from=1-2, to=1-3]
	\arrow[no head, from=1-2, to=1-3]
	\arrow[shift left=1, no head, from=2-1, to=3-1]
	\arrow[no head, from=2-1, to=3-1]
\end{tikzcd}\]
It follows from the assumptions that $F'=G''$ and $F''=G'$ (cf.~Corollary 5.7 of \cite{Batanin_Markl:kodu2022}), and we
have 
\begin{center}
\begin{tabular}{ c c c c}
 $[*_{\psi'}] \in \colorop Q(S\,F'';P')$,&
$[*_{\psi''}] \in \colorop Q(S\,F';P''),$ &
$[*_{\phi'}] \in \colorop Q(P'\,F';T),$&
$[*_{\phi''}] \in \colorop Q(P''\,F'';T).$
\end{tabular}
\end{center}
The composite $[*_{\phi'}]\circ_2[*_{\psi'}]$ is an element of~$F(Q)\colorop (S\,F''\,F';T),$
$[*_{\phi''}]\circ_2[*_{\psi''}]$
is an element of~$F(Q)\colorop (S\,F'\,F'';T),$
and hence $\tau([*_{\phi''}]\circ_2[*_{\psi''}])$
is an element of $F(Q)\colorop (S\,F''\,F';T).$
The~parallel associators are
\begin{equation}\label{equation:ideal2}
    [*_{\phi'}]\circ_2[*_{\psi'}]-\tau([*_{\phi''}]\circ_2[*_{\psi''}]),
\end{equation}
for any $\phi', \phi'', \psi'$ and $\psi''$ as above. The symbol $\tau$ stands for the isomorphism which swaps the second and third input of 
$F(Q)\colorop (S\,F'\,F'';T)$. 

Finally, we define the binary quadratic $\bO_{\it{iso}}^{\it{op}}$-operad $$\mathbb{H}_\bO:= F(Q)/\textit{As},$$ where $\textit{As}$ is the operadic ideal generated by elements (\ref{equation:ideal1}) and (\ref{equation:ideal2}).
\begin{theoremB2}\label{theorem:O-operads_as_algebras}
	The category of algebras of the binary quadratic $\bO_{\it{iso}}^{\it{op}}$-operad $\bH_{\bO}$ is isomorphic to the category of (non-unital) Markl $\bO$-operads in $\Vect$.
\end{theoremB2}
\pf
	The proof is analogous to the proof of Theorem B1. An algebra of $\bH$ is a functor 
 \[
\begin{tikzcd}
	\bO_{\it{iso}}^{\it{op}} & \Vect,
	\arrow["M", from=1-1, to=1-2]
\end{tikzcd}
\] 
together with structure maps 
 \[
\begin{tikzcd}
	M(S)\otimes M(F) & M(T)
	\arrow["\circ_{\phi}", from=1-1, to=1-2]
\end{tikzcd}
\] for any elementary morphism $\phi\colon F \triangleright T \rightarrow S$. Due to the identifications (\ref{eq:equivariance_of_*_operadic}), (\ref{equation:ideal1}), and (\ref{equation:ideal2}), the maps $\circ_{\phi}$ are equivariant and associative in the sense of (44),(46), and (47) of~\cite{Batanin_Markl:kodu2022}.
\epf
Unital Markl $\bO$-operads can be described by the operad $\bH_{\bO}$, extended by generators $$\eta_U \in \colorop \bH(\emptyset;U)$$ for every chosen local terminal object $U$ of $\bO$, with suitable relations corresponding to~(48) of~\cite{Batanin_Markl:kodu2022}.
\begin{remark}
The Theorem B1 is a direct consequence of the Theorem B2 for $\bO=\textit{Fin}$, the skeletal category of finite sets.
In \cite{Batanin_Markl:kodu2022}, it is shown that operadic categories of various kinds of graphs, such as rooted trees or genus graded connected graphs, satisfy the assumptions given in the beginning of this subsection. In this case, Markl \hbox{$\bO$-operads} are graph-indexed operads. All sorts of operad-like structures and their strongly homotopy versions can be encoded as algebras of Markl \hbox{$\bO$-operads}, cf.~\cite{Batanin_Markl_Obrad:models}. The Theorem~B2 says, that the encoding hyperoperads are themselves algebras of a binary quadratic \hbox{$\bO_{\it{iso}}^{\it{op}}$-operad}~$\bH_{\bO}$.
The Theorem B3 below is a separate result, since there is no operadic category in the classical sense (i.e.~with chosen local terminal objects for each connected component), whose operads are $C$-operads.
\end{remark}
 \subsection*{$C$-operads}
  The key ingredient in the construction of the hyperoperad $\bH_C$ is the description of $C$-operads in terms of the operations~$\circ_i^f$ of Proposition \ref{proposition:alternative_C_operad}.
The category of colors for the hyperoperad $\bH_C$ is $$Bq^{\Sigma}(C):=\bigoplus_{n\in \NN}(C^{\it{op}})^{\times n}\times C \times \Sigma_n.$$ The objects of $Bq^{\Sigma}(C)$ are 'bouquets' 
$$\colorop ( c_1 \ \cdots\ c_n  \ \sigma;c).$$
The hyperoperad $\bH_C$ will be generated by symbols $$*^f_i \in \coloropbig \bH_C({\colorop( c_1 \ \cdots\ c_n  \ \uu;c)}\ {\colorop(d_1 \ \cdots\ d_m  \ \uu;d)};{\colorop(c_1  \cdots d_1 \cdots d_m\cdots c_n  \ \uu;c)}),$$
with $f\colon d\to c_i$, equipped with free action of $(C^{\it{op}})^{\times k}\times C \times \Sigma_k$ on inputs and output, for $k=n,m$ and $n+m-1$, respectively. By the free actions we mean the left Kan extension along the discrete subcategory inclusion functor, as in~(\ref{equation:Kan_Extension}). We make the symbols \hbox{$C$-equivariant} by dividing the generating collection by the relation~(\ref{equation:C-equivariance}) of Proposition~\ref{proposition:alternative_C_operad}. The hyperoperad $\bH_C$ is obtained as a free operad, divided by the operadic ideal of \hbox{$C$-associative} relations~(\ref{equation:C-associativity}).
\begin{theoremB3}\label{theorem:C-operads_as_algebras}
For a fixed small category $C$, the category of algebras of the binary quadratic $Bq^{\Sigma}(C)$-operad $\bH_C$, described above, is isomorphic to the category of (non-unital) $C$-operads.
\end{theoremB3}	

 \begin{references*}

\bibitem[Batanin08]{Batanin:eckmann}
M.A.~Batanin.
\newblock The Eckmann--Hilton argument and higher operads.
\newblock {\em Advances in Mathematics}, 217(1):334--385, 2008.

\bibitem[BB17]{Batanin_Berger:polynomial}
M.A.~Batanin and C.~Berger.
\newblock Homotopy theory for algebras over polynomial monads.
\newblock {\em Theory and Applications of Categories}, 32(6):148--253, 2017.

\bibitem[BM15]{Batanin_Markl:duoidal}
M.A.~Batanin and M.~Markl.
\newblock Operadic categories and duoidal Deligne's conjecture.
\newblock {\em Advances in Mathematics}, 285:1630--1687, 2015.

\bibitem[BM23a]{Batanin_Markl:kodu2022}
M.A.~Batanin and M.~Markl.
\newblock Operadic categories as a natural environment for {K}oszul duality,
\newblock {\em Compositionality}, 5(3), 2023.

\bibitem[BM23b]{Batanin_Markl:kodu2021}
M.A.~Batanin and M.~Markl.
\newblock Koszul duality for operadic categories,
\newblock {\em Compositionality}, 5(4), 2023.

\bibitem[BMO23]{Batanin_Markl_Obrad:models}
M.A.~Batanin, M.~Markl, and J.~Obradovi{\'c}.
\newblock Minimal models for graph-related (hyper) operads.
\newblock {\em Journal of Pure and Applied Algebra}, page 107329, 2023.

\bibitem[BW22]{Batanin_White:substitudes}
M.A.~Batanin and D.~White.
\newblock Homotopy theory of algebras of substitudes and their localisation.
\newblock {\em Transactions of the American Mathematical Society},
  375(05):3569--3640, 2022.

\bibitem[CCN22]{CCN:moduli}
D.~Calaque, R.~Campos, and J. Nuiten.
\newblock Moduli problems for operadic algebras, 
\newblock {\em Journal of the London Mathematical Society}, 2022.

\bibitem[Day70]{Day:closed}
B.~Day.
\newblock On closed categories of functors.
\newblock In {\em Reports of the Midwest Category Seminar IV}, pages 1--38.
  Springer, 1970.

\bibitem[DS03]{Day:substitution}
B.~Day and R.~Street.
\newblock Abstract substitution in enriched categories.
\newblock {\em Journal of Pure and Applied Algebra}, 179(1-2):49--63, 2003.

\bibitem[DV21]{Dehling_Vallette:symmetric}
M.~Dehling and B.~Vallette.
\newblock Symmetric homotopy theory for operads.
\newblock {\em Algebraic \& Geometric Topology}, 21(4):1595--1660, 2021.

\bibitem[DSVV20]{DSVV}
V.~Dotsenko, S.~Shadrin, A.~Vaintrob, and B.~Vallette.
\newblock Deformation theory of cohomological field theories,
Preprint {\tt arXiv:2006.01649}, 2020.

\bibitem[FGHW08]{FGHW:Species}
M.~Fiore, N.~Gambino, M.~Hyland, and G.~Winskel.
\newblock The cartesian closed bicategory of generalised species of structures, 
\newblock {\em Journal of the London Mathematical Society}, 77(1):203--220, 2008.

\bibitem[KW17]{Kaufmann_Ward:Fey}
R.M. Kaufmann and B.C. Ward.
\newblock Feynman categories, 
\newblock{\em Ast{\'e}risque}, 387, 2017.

\bibitem[Markl96]{markl1996models}
M.~Markl.
\newblock Models for operads.
\newblock {\em Communications in Algebra}, 24(4):1471--1500, 1996.

\bibitem[Markl08]{markl2008operads}
M.~Markl.
\newblock Operads and props.
\newblock {\em Handbook of algebra}, 5:87--140, 2008.

\bibitem[MSS02]{markl2002operads}
M.~Markl, S.~Shnider, and J.~Stasheff.
\newblock Operads in algebra, topology and physics.
\newblock {\em Mathematical surveys and monographs}, 96, 2002.

\bibitem[Petersen13]{petersen2013operad}
D.~Petersen.
\newblock The operad structure of admissible G-covers.
\newblock {\em Algebra \& Number Theory}, 7(8):1953--1975, 2013.

\bibitem[Stoeckl23]{stoeckl2023koszul}
K.~Stoeckl.
\newblock Koszul Operads Governing Props and Wheeled Props,
Preprint {\tt arXiv:2308.08718}, 2023.

\bibitem[vdLaan03]{vdLaan:coloured}
P.~van~der Laan.
\newblock Coloured Koszul duality and strongly homotopy operads,
Preprint {\tt arXiv:math/0312147}, 2003.

\bibitem[Ward22]{Ward:Massey}
B.~C.~Ward.
\newblock Massey Products for Graph Homology, 
\newblock {\em International Mathematics Research Notices}, 2022(11):8086-8161, 2022.
\end{references*}
\end{document}